\documentclass[eqno,12pt]{article}

\usepackage{amssymb}
\usepackage{euscript}
\usepackage[dvips]{graphicx}
\usepackage{flafter}
\usepackage{pstricks}
\usepackage[latin1]{inputenc}
\usepackage{amsmath}
\usepackage{amsfonts}
\usepackage{pstricks-add}

\usepackage{mathrsfs} 

\evensidemargin\oddsidemargin

\begin{document}

\renewcommand{\theequation}{\thesection.\arabic{equation}}
\newtheorem{theorem}{Theorem}[section]
\newtheorem{lemma}[theorem]{Lemma}
\newtheorem{proposition}[theorem]{Proposition}
\newtheorem{corollary}[theorem]{Corollary}
\newtheorem{remark}[theorem]{Remark}
\newtheorem{fact}[theorem]{Fact}
\newtheorem{problem}[theorem]{Problem}
\newtheorem{example}[theorem]{Example}

\newcommand{\eqnsection}{
\renewcommand{\theequation}{\thesection.\arabic{equation}}
    \makeatletter
    \csname  @addtoreset\endcsname{equation}{section}
    \makeatother}
\eqnsection

\def\r{{\mathbb R}}
\def\e{{\mathbb E}}
\def\p{{\mathbb P}}
\def\q{{\mathbb Q}}
\def\qa{{\mathbb Q}^{(\alpha)}}
\def\z{{\mathbb Z}}
\def\n{{\mathbb N}}
\def\T{{\mathbb T}}
\def\ee{\mathrm{e}}
\def\L{{\pounds}}
\def\S{{\mathcal S}}
\def\F{{\mathcal F}}
\def\G{{\mathcal G}}
\def\M{{\mathbb {M}}}
\def\w{{\tt  w}}
\def\law{{\buildrel \mbox{\tiny  law} \over =}}


\centerline{ \Large \bf  How big is   the minimum of a branching random walk? }

\bigskip

\medskip

 \centerline{Yueyun Hu\footnote{Département de Mathématiques (LAGA, CNRS UMR 7539),  Université Paris XIII, Sorbonne Paris Cité, 99 avenue J.B. Clément, 93430  Villetaneuse.   Research partially  supported by ANR  (MEMEMOII)  2010 BLAN 0125 Email: yueyun@math.univ-paris13.fr}}

\medskip

 \centerline{\it Universit\'e Paris XIII}

\bigskip

{\leftskip=2truecm
\rightskip=2truecm
\baselineskip=15pt
\small

\noindent{\slshape\bfseries Summary.} Let $\M_n$ be the minimal position in the $n$-th generation, of a real-valued branching random walk in the boundary case.  As $n \to \infty$,  $\M_n- {3 \over 2} \log n$  is tight (see \cite{AbR, BZ06, A11}).  We establish here a law of iterated logarithm for the upper limits of $\M_n$: upon the system's non-extinction, $ \limsup_{n\to \infty} {1\over \log \log \log n} ( \M_n - {3\over2} \log n) = 1$ almost surely. We also  study  the  problem of  moderate deviations      of  $\M_n$: $\p(\M_n- {3 \over 2} \log n > \lambda)$ for $\lambda\to \infty$ and $\lambda=o(\log n)$. This problem is closely related to  the small deviations   of  a class of  Mandelbrot's cascades.    
\medskip

\noindent{\slshape\bfseries Résumé.}  Soit $\M_n$ la position minimale à la $n^{\mathrm {ieme}}$ génération,  d'une marche aléatoire branchante réelle dans le cas frontière. Quand $n \to \infty$,    $\M_n- {3 \over 2} \log n$  est tendue (voir \cite{AbR, BZ06, A11}).   Nous établissons une loi du logarithme itéré  pour décrire les limites supérieures de  $\M_n$: sur l'événement de la survie du système,  $ \limsup_{n\to \infty} {1\over \log \log \log n} ( \M_n - {3\over2} \log n) = 1$ presque sûrement.  Nous étudiions également les déviations modérées de  $\M_n$ : $\p(\M_n- {3 \over 2} \log n > \lambda)$ pour  $\lambda\to \infty$ et $\lambda=o(\log n)$. Ce problème est directement lié aux  petites déviations   d'une classe des cascades de    Mandelbrot.    
\medskip

  \medskip


\medskip

} 

\medskip

\section{Introduction}

Let  $\{ V(u), u \in \T\}$ be   a discrete-time  branching random walk (BRW) on the real line $\r$ driven by a point process $\Theta$.   At generation $0$, there is a single particle at the origin from which we generate a point process $\Theta$ on $\r$. The particles in $\Theta$ together with  their positions in  $\r$ constitute the first generation of the BRW.  From the position of  each particle at the first generation, we generate an independent copy of $\Theta$. The collection of all particles together with their positions gives  the second generation of the BRW, and so on.    The genealogy of all particles forms  a Galton-Watson tree $\T$ (whose root is denoted by $\varnothing$). For any particle $u \in \T$, we denote by $V(u )$ its position in  $\r$ and $\vert u \vert $ its generation in $\T$.  The whole system may die out or survive forever.

Plainly    $\Theta= \sum_{ \vert u \vert =1} \delta_{ \{ V(u)\}}$.   Let $\nu= \Theta(\r)$.   Throughout this paper and unless stated otherwise, we shall assume  that the BRW is  in the boundary case, i.e.   \begin{equation}\label{hyp}  \e  [\nu]  \in (1, \infty], \qquad \e \Big[  \sum_{ \vert u \vert =1} \ee^{- V(u)} \Big]=1 , \qquad \e \Big[  \sum_{ \vert u \vert =1}  V(u) \, \ee^{- V(u)} \Big] =0. \end{equation}

\noindent Notice  that  under \eqref{hyp}, it is possible   that    $\p(\nu=\infty ) >0$.  See  Jaffuel \cite{J12} for detailed discussions on how to reduce   a general branching random walk   to the boundary case.

Denote by $\M_n:= \min_{ \vert u \vert = n} V(u)$  the minimum  of the branching random walk in the $n$-th generation (with convention: $\inf\emptyset  \equiv \infty$).   Hammersly \cite{H74}, Kingman \cite{K75} and Biggins \cite{B76} established the law of large numbers for $\M_n$ (for any general branching random walk), whereas the second order limits have   attracted many recent  attentions, see \cite{AbR, HS09, BZ06, A11} and the references therein.  In particular,  A{\"{\i}}d{\'e}kon \cite{A11}     proved the convergence in law of $\M_n- {3\over 2} \log n$ under \eqref{hyp} and some mild conditions.

On   the almost sure limits of $\M_n$, it was shown in   \cite{HS09} that  there is  the following phenomena of fluctuation at the logarithmic scale.  Assume \eqref{hyp}.  If there exists some  $\delta>0$ such that  $\e[ \nu^{1+\delta}] < \infty$ and $\e \big[ \int_{\r} (\ee^{ \delta x } + \ee^{- (1+\delta)x } ) \Theta(dx)\big] < \infty$,  then  \begin{eqnarray*} \limsup_{ n \to \infty} { \M_n \over \log n}  =  { 3\over 2}   \qquad \ \mbox{and } \qquad 
  \liminf_{ n \to \infty} { \M_n \over \log n}  =  { 1\over 2}, \qquad \p^*\mbox{-a.s.}, 
\end{eqnarray*}

\noindent where here and in the sequel,      $$ \p^*(\cdot):= \p \left( \cdot \vert \S \right),$$ and $\S:=\{\T \mbox{ is not finite}\} $ denotes the event that the whole system survives.

It turns out that much more can be said on the  lower  limits  ${1\over2} \log n  $ of   $\M_n$:   Under   \eqref{hyp} and the following integrability condition \begin{equation} \label{int1}
     \sigma^2:=   \e \Big[  \sum_{ \vert u \vert =1} (V(u))^2\, \ee^{- V(u)} \Big] < \infty, \quad  \e \Big[ \zeta ( (\log \zeta)^+)^2 + \widetilde \zeta (\log \widetilde \zeta)^+  \Big] < \infty,  
\end{equation} with    $\zeta:=   \sum_{ \vert u \vert =1} \ee^{- V(u)}  $, $ \widetilde \zeta:=  \sum_{ \vert u \vert =1}  (V(u))^+\, \,  \ee^{- V(u)}   $   and $x^+:= \max( 0,  x)$,   A{\"{\i}}d{\'e}kon and Shi \cite{AS12} proved that   $$ \liminf_{ n\to \infty} \Big( \M_n - { 1\over 2} \log n \Big) = - \infty, \qquad \mbox{$\p^*$-a.s.}$$
 
 \noindent 
 Furthermore, by following A{\"{\i}}d{\'e}kon and Shi \cite{AS12}'s methods, we     established (\cite{H12}) an integral test       to describe the lower limits of $\M_n- {1\over 2} \log n$. As a consequence,   we have that under \eqref{hyp}  and  \eqref{int1},   \begin{equation}\label{lowerLIL} \liminf_{ n \to \infty} { 1\over  \log \log n}  \Big( \M_n - { 1\over 2 }\log n \Big) =-1, \qquad \mbox{$\p^*$-a.s.}\end{equation}

 In this paper, we wish to investigate  how big  $\M_n - {3\over 2} \log n$ can be.  The following law of iterated logarithm (LIL) describes the upper limits    of $\M_n$: 
 
  \begin{theorem}  \label{t:3log}   Assume \eqref{hyp}, \eqref{int1} and  that $\e \big[ \sum_{\vert u \vert=1} (V(u)^+)^3 \ee^{- V(u)}\big] < \infty$.  We suppose that the law of $\log \sum_{|u|=1} \ee^{-V(u)}$ is non-lattice.  Then  \begin{equation} \limsup_{n\to \infty} {1\over \log \log \log n} ( \M_n - {3\over2} \log n) = 1, \qquad \p^* \mbox{-a.s.} \end{equation}        \end{theorem}

  The  integrability of  $ \sum_{\vert u \vert=1} (V(u)^+)^3 \ee^{- V(u)} $  is   used only in the proof of Lemma \ref{L:Nerman},  see Remark \ref{remark-nerman}, Section \ref{S:4}. 
  
  The assumption of the non-lattice law of $\log \sum_{|u|=1} \ee^{-V(u)}$ is used  only in the proof of Lemma \ref{tightness}, see the footnote therein.

 Usually,  to establish such LIL, the first step would be the study of the moderate deviations: $$ \p^* \Big( \M_n - {3 \over 2} \log  n > \lambda \Big), \qquad   \mbox{when } \lambda= o(\log n) \mbox{ and } \lambda, n   \to \infty.$$

  
      Denote by  $p_j= \p(\nu=j), j \ge0$,  the offspring distribution of the Galton-Watson tree $\T$.   Concerning the small deviations of the size of $\T$, there exist  two cases: either $p_0+p_1>0$ (namely  the Schr\"{o}der case)  or $p_0=p_1=0$  (namely  the B\"{o}ttcher case), see e.g.   Fleischmann and Wachtel \cite{FW07, FW09} and the references therein.   Basically in the Schr\"{o}der case, the  tree $\T$ may grow linearly whereas it always grows exponentially  in the B\"{o}ttcher case.    For the branching random walk, we        shall prove  that   the moderate deviations of $\M_n$  decay  exponentially fast  or double-exponentially fast  depending on  the growth rate of  $\T$.

   Let $q:= \p(\T \mbox{ is finite}) = \p( \S^c)\in [0, 1)$ be the extinction probability. We introduce two separate cases:  {\sl 
 \begin{description}
\item[(the Schr\"{o}der case)] if the following hypotheses hold:  \begin{equation} \label{defgamma1} \e \Big [ 1_{(\nu\ge 1)}\,  q^{ \nu-1} \,   \sum_{ \vert u \vert =1} \ee^{\gamma\, V(u)}   \Big]  = 1, \qquad \mbox{for some constant $\gamma >0$,} \end{equation}
\noindent and  \begin{equation} 
  \label{hyp-exp} \e \Big[ \sum_{ \vert u \vert =1} \ee^{ a V(u)} \Big]   <  \infty, \qquad \mbox{for some $a>\gamma$}.   \end{equation}     
    \item[(the B\"{o}ttcher case)]  if the following hypotheses hold:  \begin{eqnarray}\label{bot}    &&p_0  = p_1=0, \\
  \label{bounded} &&   \sup_{ \vert u \vert =1} V(u) \le K, \qquad \mbox{for some constant $K>0$.} \end{eqnarray}    \end{description}
}

\begin{remark}\label{R:1}     
$\phantom{vide}$

(i)  When a.s. $\nu\ge1$   in  the Schr\"{o}der case,   the condition \eqref{defgamma1} just amounts to    \begin{equation}
\e \Big[ 1_{(\nu=1)} \, \sum_{ \vert u \vert =1} \ee^{\gamma\, V(u)}   \Big]  =  1, \qquad \mbox{if } q =0. \label{defgamma2}\end{equation}

(ii) Under \eqref{hyp},  the condition   \eqref{hyp-exp} or \eqref{bounded} implies  that $\e[\nu] <\infty$.  The technical conditions \eqref{hyp-exp} and \eqref{bounded} are made   to  avoid  too large jumps of $\Theta$ in the   moderate deviations. 	

(iii)   In the  B\"{o}ttcher case,  we  can define   a parameter $\beta>0$ by   \begin{equation}\label{defbeta}   \beta:= \sup\{ a>0: \p\big( \sum_{ \vert u \vert =1} \ee^{- a \,  V(u)} \ge 1\big)=1\}.
   \end{equation} 
\noindent  Note that $\beta<1$ if we assume \eqref{hyp}.    \end{remark}



 The parameters $\gamma$ and $\beta$     will   naturally  appear in  the  small deviations of a class of Mandelbrot's cascades.  Under \eqref{hyp} and \eqref{int1}, the so-called derivative martingale (with  convention: $\sum_\emptyset :=0$) $$D_n:= \sum_{ \vert u \vert = n} V(u) \, \ee^{- V(u)}, \qquad n \ge 0,$$ 
 
 \noindent  converges almost surely to some limit $D_\infty$ which is $\p^*$-a.s. positive (see e.g. Biggins and Kyprianou \cite{BK04} and  A{\"{\i}}d{\'e}kon \cite{A11}).   The nonnegative random  variable $D_\infty$ satisfies the following equation in law (Mandelbrot's cascade):  \begin{equation}\label{eq-cascade} D_\infty \, \law\,  \sum_{\vert u \vert =1} \ee^{- V(u)} D^{(u)}_\infty,\end{equation}

\noindent where conditioned on $\{ V(u), \vert u \vert =1\}$,  $(D^{(u)}_\infty)_{\vert u \vert =1}$ are independent copies of $D_\infty$. The moderate deviations of $\M_n$ will be naturally related to the small deviations   of $ D_\infty  $ which were  already studied in the literature, see e.g. Liu \cite{Liu99, Liu01} and the references therein.

We shall work under a more general setting in order    that  Theorem \ref{T:smalldev} could  also be applied to  the non-degenerated case  of Mandelbrot's cascades.  Instead of \eqref{hyp}, we assume that there exists some constant $\chi \in (0, 1]$ such that \begin{equation}\label{cascade-hyp} \e \Big[ \sum_{\vert u \vert =1} \ee^{- \chi V(u)} \big] \le 1,   \qquad  \mbox{and} \qquad    \e [\nu] \in (1, \infty],
\end{equation}

 \noindent where as before, $\nu:= \sum_{\vert u \vert =1} 1$. 
 
 The condition \eqref{cascade-hyp} ensures that there exists a non-trivial nonnegative solution $Z$ to the following equation: \begin{equation}\label{cascade1} Z \, \law\, \sum_{\vert u \vert =1} \ee^{- V(u)} Z^{(u)}, 
\end{equation} where conditioned on $\{ V(u), \vert u \vert =1\}$,  $(Z^{(u)})_{\vert u \vert =1}$ are independent copies of $Z$, see  Liu \cite{Liu01}, Proposition 1.1.

Denote by $f(x) \asymp g(x)$ [resp: $f(x) \sim g(x)$] as $x \to x_0$ if $0 < \liminf_{x \to x_0} f(x)/g(x) \le \limsup_{x\to x_0} f(x)/g(x) < \infty$ [resp: $\lim_{x\to x_0} f(x)/g(x)=1$].  The following result may arise  an   interest in  Mandelbrot's cascades.

\begin{theorem}\label{T:smalldev}  Assume \eqref{cascade-hyp}. Let $Z\ge0$ be a  non-trivial solution of \eqref{cascade1}. 

\begin{description}
\item[(The Schr\"{o}der case)]  Assume   \eqref{defgamma1}   and     \eqref{hyp-exp}. Then \begin{equation}\label{smalldev1new}    \p\Big(  0< Z < \varepsilon \Big)  \, \asymp\, \varepsilon^\gamma, \qquad  \mbox{as }  \varepsilon \to 0,      \end{equation}  and $\e \big[ \ee^{- t Z } 1_{( Z>0)}\big] \asymp  t^{- \gamma  }$ as $t \to \infty$.  
\item[(The B\"{o}ttcher case)] Assume   \eqref{bot},     \eqref{bounded} and that $\sum_{\vert u \vert=1} \ee^{- \chi V(u)} \not \equiv 1$. Then  \begin{equation} \label{smalldev2}  \e \Big[ \ee^{- t Z } \Big]= \ee^{- t^{ \beta+ o(1) }}, \quad t \to \infty, \end{equation} and $\p\big( Z  < \varepsilon\big) = \ee^{ - \varepsilon^{- {\beta \over 1-\beta} + o(1)}}, $  as $ \varepsilon \to 0$, with $\beta$ defined in \eqref{defbeta}.
\end{description}

\end{theorem}

 
Obviously  we can apply Theorem \ref{T:smalldev} to  $Z:=D_\infty$ with $\chi=1$.   In the B\"{o}ttcher case, the two conditions \eqref{cascade-hyp} and $\sum_{\vert u \vert=1} \ee^{- \chi V(u)} \not \equiv 1$  imply  that  $\beta < \chi $, hence   $\beta < 1$; moreover, $\mbox{essinf}  \, \sum_{ \vert u \vert =1} \ee^{- \beta V(u)} = 1$.   

 Let us mention that  \eqref{smalldev1new} confirms a prediction  in Liu \cite{Liu01} who already proved that if $q=0$, then for any $a>0$,  $\e \big[  Z^{-a}\big] < \infty$ if and only if $a < \gamma$. When all $V(u), \vert u \vert =1$,  are equal  to some random variable,  \eqref{smalldev2} is in agreement with Liu \cite{Liu99}, Theorem 6.1. If  furthermore, all $V(u)$ are equal to some constant,   then \eqref{smalldev1new} and \eqref{smalldev2} give   some rough estimates on the   limiting  law of   Galton-Watson processes, see Fleischmann and Wachtel \cite{FW07}, \cite{FW09}   for  the precise  estimates.  We refer to    \cite{BGMS} for further studies of the conditioned Galton-Watson tree   itself.     For instance, we could   seek  the asymptotic behaviors of the BRW conditioned on $\{0<D_\infty<\varepsilon\}$, as $\varepsilon \to 0$, but this problem  exceeds the scope of the present paper.

Our moderate deviations result on $\M_n$ reads  as follows:

    \begin{theorem} \label{T:DEV1} Assume \eqref{hyp}, \eqref{int1}. Let $\lambda, n \to \infty$ and $\lambda = o( \log n)$. 
    \begin{description}
\item[(The Schr\"{o}der case)]  Assume    \eqref{defgamma1}  and that \eqref{hyp-exp} hold for all $a>0$. Then \begin{equation}\label{dev1-main} \p^* \Big(   \M_n > { 3\over2} \log n +   \lambda \Big) = \ee^{- (\gamma+ o(1)) \lambda}. \end{equation}    \item[(The B\"{o}ttcher case)] Assume   \eqref{bot} and     \eqref{bounded}. Then \begin{equation}\label{dev2-main} \p\Big(   \M_n > { 3\over2} \log n +   \lambda \Big) = \exp( - \ee^{(\beta+ o(1)) \lambda}). \end{equation}  \end{description}  The same estimates hold  if we replace $\M_n$ by $ \max_{ n \le k \le 2 n}  \M_k$.
   \end{theorem}

We refer to  A{\"{\i}}d{\'e}kon \cite{A11}, Proposition 4.1 for the precise estimate on $\p(\M_n < {3\over2} \log n - \lambda)$ as $ \lambda \le {3\over2}    \log n$ and $\lambda \to \infty$. 

  Comparing   Theorem \ref{t:3log} and  Theorem  \ref{T:DEV1}, we remark that the almost sure behaviors of $\M_n$ are not related to the moderate deviations of $\M_n$.  This can be explained as follows: Define for all $ \lambda\ge0$ and $u \in \T$,  \begin{equation} \label{tau} \tau_\lambda(u):= \inf\{ 1 \le i \le \vert u \vert : V(u_i) > \lambda\},  \qquad  \mbox{(with convention $ \inf \emptyset= \infty$)}, 
\end{equation} 

\noindent  where here and in the sequel, $\{u_0= \varnothing, u_1, ..., u_{\vert u \vert}:=u\}$ denotes the shortest path from $\varnothing$ to $u$ such that $\vert u_i\vert =i$ for all $0 \le i \le \vert u \vert$. We introduce the stopping lines:  \begin{equation}  \label{defll} \L_\lambda  :=   \{ u \in \T : \tau_\lambda(u)= \vert u \vert\}, \qquad \lambda\ge0.    \end{equation}

Roughly speaking, the almost sure   limits of $\M_n$ ($\limsup$ of $\M_n$)  are determined by those of $\#\L_\lambda$, whereas the moderate deviations of $\M_n$ are by the small deviations   of $\#\L_\lambda$.  By Nerman \cite{N81}, $\p^*$-almost surely, $\#\L_\lambda$ is of order $\ee^{ (1+o(1) )\lambda}$; however, to make $\#\L_\lambda$ to be as small as possible (and conditioned on $\{ \#\L_\lambda >0\}$),  in the Schr\"{o}der case,  $\L_\lambda$ will be essentially a singleton or a set of few points  with exponential costs (see Lemma \ref{lowerproba=1}), which is no longer   possible in the B\"{o}ttcher case. To relate $\#\L_\lambda$ to $D_\infty$, we shall use the martingale $(D_n)$ at the stopping line $\L_\lambda$: \begin{equation}\label{DL} D_{\L_\lambda}:= \sum_{ u \in \L_\lambda} V(u) \ee^{- V(u)}, \end{equation}

\noindent which, as shown in Biggins and Kyprianou \cite{BK04}, converges almost surely to $D_\infty$ as $\lambda \to \infty$.   For $u \in \L_\lambda$,   $V(u) \approx \lambda$, hence  $D_{\L_\lambda} \approx  \lambda \, \ee^{-\lambda} \#\L_\lambda$.  Then the problem of small values of $\#\L_\lambda$ will be reduced to that of $D_{\L_\lambda}$ and $D_\infty$ as $\lambda \to \infty$. The hypothesis \eqref{hyp-exp} and \eqref{bounded} are made  to control the possible overshoots.

  The rest of the paper is organized as follows:    In Section \ref{S:2}, we collect some facts on a one-dimensional random walk and on the branching random walk.  In Section \ref{S:3}, we study  the   cascade equation \eqref{cascade1} and  prove Theorem \ref{T:smalldev}.    In Section \ref{S:4}, we first prove some  uniform   tightness of $\M_n-{3\over2}\log n$  (Lemma \ref{tightness}) and  then Theorem \ref{t:3log}. Finally, in Section \ref{S:5}, we prove Theorem \ref{T:DEV1} in two separate subsections on  the Schr\"{o}der case and on the B\"{o}ttcher case.

  Throughout the paper,  we adopt the usual conventions that $\sum_\emptyset:=0$, $\sup_\emptyset:= 0$, $\prod_\emptyset:=1$, $\inf_\emptyset:=\infty$; we also denote by  $(c_i, 1\le i \le 15)$    some positive constants,  and  by  $C, C'$ and $C^{''}$ (eventually with a subscript)    some unimportant positive constants whose values can vary from one paragraph to another one.

\section{Preliminaries} \label{S:2}

\subsection{Estimates on a centered real-valued random walk}

We collect here some estimates on a   real-valued random walk $\{S_k, k\ge0\}$, under $\p$, centered and with finite   variance $\sigma^2 >0$. Write $\p_x$ and $\e_x$  when $S_0=x$.  Let $\underline S_n:=  \min_{  0\le i \le n} S_i$, $\forall\, n \ge 0$. The renewal function $R(x)$ related to the random walk $S$  is defined as follows: \begin{equation}\label{rx}  R(x):= \sum_{k=0}^\infty \p \Big( S_k \ge -x , \, S_k < \underline S_{k-1} \Big) , \qquad x \ge0, \end{equation} and $R(x)=0$ if $x <0$.    Moreover (see Feller \cite{feller}, pp.612),  \begin{equation}\label{defcr} \lim_{x \to \infty} { R(x)\over x}= c_1>0 . \end{equation}

\begin{lemma} \label{L:fact}   Let $S$ be a centered random walk with finite and positive  variance. There exists some constant $c_2>0$  such that for any $b \ge a \ge 0,  x\ge0,  n\ge 1$,   \begin{equation} \label{AS1} \p_x \Big( S_n \in [a, b],  \underline S_n  \ge 0\Big) \le c_2  \,  (1+x)  (1+b- a) (1+b) n^{-{3\over2}} . \end{equation}

\noindent
For any fixed $0< r< 1$, there exists some $c_3\equiv c_{3, r} >0$ such that for all  $b \ge a \ge 0$, $x,   y \ge 0$ and $ n\ge 1$, \begin{eqnarray}  \label{AS2}  \p_x \Big( S_n \in [y+a, y+b],  \underline S_n  \ge 0, \min_{ r   n \le i <  n} S_i \ge y \Big)  &\le & c_3  \, (1+x)  (1+b-a) (1+b) n^{-{3\over2}}, \\
    \label{ladder1} 
  \p_x \Big( \underline S_n \ge 0,  \min_{ r n  \le i < n} S_i > y,  S_n \le y \Big) & \le& c_3 \, (1+x)   n^{-{3\over2}}. \end{eqnarray}
\noindent
For any $a >0$, if $\e\big[ S_1 ^2 \ee^{ a S_1}\big] < \infty$, then there exists some $C_a>0$ such that for any $b \ge 0$, \begin{equation} \label{overshoot} \p \Big( S_{ \tau_b} - b > x \Big) \le C_a \, \ee^{ - a x}, \qquad \forall\, x \ge 0, \end{equation}  where $\tau_b:= \inf\{ j\ge 0: S_j > b\}$.
\end{lemma}

{\noindent \bf Proof of Lemma \ref{L:fact}.} See  A{\"{\i}}d{\'e}kon and Shi \cite{AS12} for \eqref{AS1} and  \eqref{AS2}.   To get  \eqref{overshoot},  note  that $\e\big[ S_1  ^2 \ee^{ a S_1 }\big] < \infty $ if and only if $\e\big[ (S_1^+)^2 \ee^{ a S_1^+}\big] < \infty$. By Doney (\cite{D80}, pp.250),  this condition ensures that $\e\big[ S_{\tau_0} \ee^{ a S_{\tau_0}}\big]< \infty$. Then  in view of Chang (\cite{Chang94},  Proposition 4.2), we have  that uniformly in $b>0$,  $ \e \big[ \ee^{ a (S_{\tau_b} - b)} \big] \le C_a$ for some constant $C_a>0$, which implies   \eqref{overshoot} by Chebychev's inequality.    

It remains to check \eqref{ladder1}.  Let $f(x):= \p( S_1 \le -x), x\ge0$.  It follows from the Markov property at $n-1$ that  the probability in LHS of \eqref{ladder1} equals \begin{eqnarray*}  && \e_x \Big[ 1_{( \underline S_{n-1} \ge 0, \min_{ r n \le i < n} S_i >y)} \,  f(S_{n-1} - y)    \Big] \\
     &\le&  \sum_{j=0}^\infty \,  f(j)\,  \p_x\Big( \underline S_{n-1} \ge 0, \min_{ r n \le i \le n-1} S_i >y, y + j < S_{n-1} \le y+j +1\Big) \\
     &\le& C \, (1+x) \, n^{-3/2} \, \sum_{j=0}^\infty  \, f(j)  \, (2+j)  \qquad \mbox{ (by \eqref{AS2})}\\
     &\le&  C' \, (1+x) \, n^{3/2},
\end{eqnarray*} yielding \eqref{ladder1}. $\Box$

\subsection{Change of measures   for the branching random walk}\label{S:2.2}

In this subsection, we recall some change  of measure formulas  in  the branching random walk, for the details we refer to \cite{BK04,CRW,L97, AS12,  HR, Shi} and the references therein.  

At first let us fix some notations: For $ \vert u \vert = n$, we write as before $  \{ u_0:= \varnothing, u_1, ..., u_{n-1}, u_n =u\}$ the   path  from the root $\varnothing$  to $u$ such that $\vert u_i\vert =i$ for any $0\le i \le n$.   Define $\overline V(u):= \max_{ 1\le i \le n } V(u_i)$ and $\underline V(u):= \min_{ 1 \le i \le n} V(u_i)$. For any $u, v \in \T$, we use the partial order  $ u < v$ if $u$ is an ancestor of $v$ and  $u \le v$ if $u < v$ or $u=v$.  We also denote by   ${\buildrel \leftarrow \over u}$   the parent of $u$ and by $\nu(u)$ the number of children of $u$.   Define $ \mho( u):=\{ v: {\buildrel \leftarrow \over v}= {\buildrel \leftarrow \over u}, v \neq u\}$   the set (eventually empty) of  brothers of $u$  for any $ u \not= \varnothing$.   For any $u \in \T$, we denote by $\T_u:= \{ v \in \T: u \le  v\}$ the subtree of $\T$ rooted at $u$.

Under \eqref{hyp}, there exists a centered real-valued random walk $\{S_n, n\ge0\}$ such that for any $n \ge1$ and any measurable $f:\r^n \to \r_+$ , \begin{equation} \label{many} \e \Big[ \sum_{ \vert u \vert =n} \ee^{- V(u)}  f( V(u_1), ..., V(u_n)) \Big]= \e \left( f(S_1, ..., S_n)\right), \end{equation} 

\noindent   which is often referred as  the ``many-to-one"  formula. Moreover under \eqref{int1},   $\mbox{Var}(S_1)=\sigma^2= \e \Big[ \sum_{\vert u \vert =1} (V(u))^2 \ee^{- V(u)} \Big] \in (0, \infty).$  We shall use the notation \begin{equation} \label{tau0}\tau_0:= \inf\{ j\ge 1: S_j >0\}. \end{equation}

Denote by $(\F_n, n\ge0)$ the natural filtration of the branching random walk.  Under \eqref{hyp},   the process $W_n:= \sum_{ \vert u \vert =n} \ee^{- V(u)}$, $n \ge 1$, is a $(\p , (\F_n))$-martingale.  It is well-known (see \cite{BK04,CRW, L97, AS12, HR,  Shi}) that  on some enlarged probability space (more precisely on the space of marked trees enlarged by an infinite ray $(\w_n, n \ge0)$, called spine), we may construct a probability $\q$ such that    the following statements (i), (ii) and (iii) hold:

(i)  For all $n\ge1$,  $$     { d \q  \over d \p  } \big\vert _{\F_n}=  W_n, \qquad \mbox{ and } \qquad 
    \q  \left( \w_n = u \big\vert \F_n\right)= {1\over W_n } \ee^{-V(u)}  , \qquad \forall \vert u \vert =n. $$

(ii)  Under $\q $, the process $\{V(\w_n), n \ge0\}$ along the spine $(\w_n)_{n\ge0}$,  is distributed as the random walk $(S_n, n\ge0)$ under $\p$. Moreover, $(   \sum_{ u \in \mho( \w_k)} \delta_{\{ \Delta V(u)\}},  \Delta V(\w_k))_{k \ge 1}$ are i.i.d. under $\q$, where $\Delta V(u):= V({\buildrel \leftarrow \over u})- V(u)$ for any $ u \not= \varnothing$.  

(iii) Let $\G_n:= \sigma\{ u, V(u): {\buildrel \leftarrow \over u} \in \{\w_k, 0\le k < n\}\},  $ $n\ge0$.  Then $\G_\infty$ is the $\sigma$-algebra  generated by the   spine. Under $\q $ and conditioned on $\G_\infty$, for all $u \not \in \{ \w_k, k\ge0\}$ but ${\buildrel \leftarrow \over u} \in \{\w_k,  k \ge0\}$ the induced branching random walk $(V(uv), \vert v \vert \ge0)$ are independent and are distributed as $\p_{V(u)}$, where $\{uv, \vert v \vert \ge0)$ is  the subtree  $\T_u$.

We mention   that  the above change of measure still holds  for the stopping line $\L_\lambda$ (see e.g.  \cite{EOY}, Proposition 3, for  the detailed statement): i.e.   replace $\vert u \vert=n$ by   $u \in \L_\lambda$, $\F_n$ by $\F_{\L_\lambda}$ the $\sigma$-filed generated by the BRW up to $\L_\lambda$, and   $W_n$ by  \begin{equation}\label{WL} W_{\L_\lambda}:= \sum_{u \in \L_\lambda} \ee^{- V(u)}.\end{equation}

 For brevity, we shall write    $\q[X]$  for the expectation of some random variable $X$ under the probability $\q$.

\section{Proof of Theorem \ref{T:smalldev}}\label{S:3}

The  following  result  is due to   Liu \cite{Liu01}:

    \begin{lemma}[Liu \cite{Liu01}] \label{L:liu} Assuming   \eqref{defgamma1},   \eqref{hyp-exp}  and \eqref{cascade-hyp}. Let $Z\ge 0$ be a non-trivial solution of \eqref{cascade1}.    For any $0< \varepsilon < \gamma$, there exists some positive constant $c_4=c_4(\varepsilon) $ such that   \begin{equation}\label{smalldev1}  \e \Big[ \ee^{- t Z} 1_{( Z>0)} \Big]  \le c_4\,  t^{-\gamma+\varepsilon}, \qquad \forall \,  t \ge  1.        \end{equation}
\end{lemma}
  
{ \noindent\bf Proof of Lemma \ref{L:liu}.}  At first we remark that \begin{equation}\label{probaz=0} \p \Big( Z=0\Big)= q.\end{equation}

In fact, we easily deduce from \eqref{cascade1} that  the probability $\p(Z=0)$ is a solution of  $x= \e [x^\nu]$ which only has two solutions $q$ and $1$ for  $x\in [0, 1]$. This gives \eqref{probaz=0}. 

  In the case $q=0$, namely $Z>0$ a.s.,  $\gamma$ is defined through \eqref{defgamma2}, it is easy to check that $\p(\sum_{|u|=1} \ee^{-V(u)}\not=1)>0$, then   \eqref{smalldev1}  follows  exactly  from   Liu \cite{Liu01}, Theorem 2.4, after a standard Tauberian argument (see Lemma 4.4 in \cite{Liu99}).   We only need to check that the case $q>0$ can be reduced  to the  case $q=0$.

  For   brevity,  let us  denote by $\{ A_i, 1 \le i \le \nu\}$ the family $\{ \ee^{- V(u)}, \vert u \vert =1\}$ [the order of $A_i$ is arbitrary].  Then $Z$ satisfies the equation in law  \begin{equation} \label{eq-law}  Z \, \law\,  \sum_{i=1}^\nu  A_i Z_i,\end{equation}

 \noindent  
with $(Z_i, i \ge1)$ independent copies of $Z$, and independent of $(A_i)_{1 \le i \le \nu}$.  Let $\{\xi,  \xi_i, i \ge 1\}$ be a family of i.i.d. Bernoulli random variables, independent of everything else, with common law $\p( \xi= 0)= q = 1- \p(\xi=1)$.   Let  $\widehat Z$ be a  random variable   distributed as     $Z$ conditioned on  $\{ Z > 0\}$. Since $\p(Z>0) = 1-q$, we have that  $Z \law \,  \xi \,  \widehat Z $.   Then  we deduce from \eqref{eq-law} that $$\widehat Z \, \law  \,\mbox{ $\sum_{i=1}^\nu\,  A_i \, \xi_i  \, \widehat Z_i $  conditioned on $\{ \sum_{i=1}^\nu \xi_i >0\}$ },$$

\noindent where $(\widehat Z_i, i \ge 1)$  are  i.i.d. copies of $\widehat Z$, and $(\nu, A_i, 1 \le i \le \nu)$ and $(  \xi_i, i \ge 1)$ are three independent families of random variables.   Let $\{ \widehat A_i, 1 \le i \le \widehat \nu)$ be a family of random variables such that for any nonnegative measurable function $f$, \begin{equation}\label{newai23}  \e \Big[ \ee^{ - \sum_{i=1}^{\widehat \nu}  f( \widehat A_i)} \Big]= \e \left[ \ee^{ - \sum_{i=1}^{  \nu}   \xi_i \,   f(   A_i)} \, \Big \vert   \sum_{i=1}^\nu \xi_i >0\right].\end{equation}

In other words, $\sum_{i=1}^{\widehat \nu}  \delta_{\{ \widehat A_i\}}$ has the same law as the point process $\sum_{1\le i\le   \nu, \, \xi_i\neq 0}  \delta_{\{   A_i\}}$ conditioning  the latter does not vanish everywhere.  Elementary calculations show that $\p(\sum_{i=1}^\nu \xi_i >0)= 1- \e [ q^\nu]= 1-q$ and    for   any nonnegative measurable function $f$,  \begin{eqnarray}    \e \Big[   \sum_{i=1}^{\widehat \nu}  f( \widehat A_i) \Big]  =  \e \left[   \sum_{i=1}^{  \nu}   \xi_i \,   f(   A_i) \, \big \vert   \sum_{i=1}^\nu \xi_i >0\right]  = {1\over 1-q}  \e \Big[   \sum_{i=1}^{  \nu}   \xi_i \,   f(   A_i)  \Big] = \e \Big[   \sum_{i=1}^{  \nu}   \,   f(   A_i)  \Big]. \label{newai}
\end{eqnarray}

\noindent In particular, $ \e \Big[   \sum_{i=1}^{\widehat \nu}   \widehat A_i ^\chi \Big]  = \e \Big[   \sum_{i=1}^{  \nu}     A_i ^\chi \Big]  \le 1$ and $ \e \big[  \widehat \nu \big]=\e \big[ \nu \big] \in (1, \infty]$.    Moreover,  we deduce from \eqref{newai23}   that $\widehat \nu$ is distributed as $ \sum_{i=1}^\nu \xi_i $  conditioned on $\{ \sum_{i=1}^\nu \xi_i >0\}$, hence $\widehat \nu \ge 1$ a.s.  It is easy (e.g. by using the Laplace transform) to see that $$ \widehat Z \,\law\, \sum_{i=1}^{ \widehat \nu} \widehat A_i\,  \widehat Z_i.$$  

\noindent  Therefore  we  can apply the case $q=0$ of   \eqref{smalldev1} to $\widehat Z $ once we have determined the corresponding parameter $\gamma$ (as in \eqref{defgamma2}) for $\widehat Z$. To this end,  let $t_\xi= \inf\{ 1 \le i \le \nu: \xi_i=1\}$. Then $\widehat A_1= A_{t_\xi}$ if $t_\xi < \infty$. We have      \begin{eqnarray*}  \e \Big[  (\widehat A_1)^{- \gamma} 1_{( \widehat \nu=1)}\Big] &=& \e \left [  A_{t_\xi} ^{-\gamma} 1_{ (\sum_{i=1}^\nu \xi_i =1)} \big \vert\, \sum_{i=1}^\nu \xi_i >0\right] \\
	&=& {1\over 1-q} \e \left [1_{(\nu\ge 1)}  \sum_{k=1}^\nu  A_k ^{-\gamma} 1_{( \xi_k=1, \, \xi_i=0  , \, \forall   i \not =k, 1 \le i \le \nu)} \right] \\
	&=& \e \left[ 1_{(\nu\ge 1)} q^{\nu-1} \, \sum_{k=1}^\nu \, A_k^{-\gamma}\right]  =   \e \left[ 1_{(\nu\ge 1)} q^{\nu-1} \, \sum_{\vert u \vert =1}  \, \ee^{\gamma V(u) }\right] = 1,
\end{eqnarray*}

\noindent by \eqref{defgamma1}.   Therefore    $ \e \big[ \ee^{- t \widehat Z } \big] = O(t^{-\gamma+\varepsilon})$ as $t  \to \infty$.              The Lemma follows from the fact that $ \p( 0< Z < x)= (1-q) \p( \widehat Z < x)$ for any $x>0$. $\Box$

\subsection{Proof of Theorem \ref{T:smalldev}:  the Schr\"{o}der case}

  As shown in the proof of Lemma \ref{L:liu},  we can assume $q=0$ (hence we assume \eqref{defgamma2}) in this proof without any loss of generality. 
     Let $\Phi(t):= \e \big[ \ee^{- t Z} \big]$ for $t \ge0$.  We are going to prove that  \begin{equation} \label{PhiPhi} \Phi(t) \asymp t^{- \gamma}, \qquad t \to \infty.
\end{equation}

To this end,   we have by \eqref{eq-law} that  \begin{equation}\label{Phi} \Phi(t) = \e \Big[ \prod_{i=1}^\nu \Phi(t A_i)\Big], \qquad t \ge0.\end{equation}
      
      \noindent Note also that the condition  \eqref{defgamma2} can be re-written as $ \e \big[ 1_{(\nu=1)} A_1^{- \gamma}\big]=1$.  Define $g(t):= t^\gamma \Phi(t)$ for all $t \ge0$.  Then  for any $t>0$,  \begin{equation}\label{ite1*} g(t)= t^\gamma\; \Phi(t) \ge  \, t^\gamma \e \Big[ 1_{(\nu=1)} \Phi(t A_1)\Big] = \e \Big[ 1_{(\nu=1)}  A_1^{-\gamma} \, g( t A_1)\Big] = \e \Big[ g(t \widetilde A_1)\Big],\end{equation}

     \noindent where $\widetilde A_1$ denotes a (positive) random variable whose law is determined by $ \e \big[ f( \widetilde A_1)\big] :=  \e \big[ 1_{(\nu=1)}  A_1^{-\gamma} \, f(  A_1)\big] $ for any measurable bounded function $f$.  In particular,  $ \e \big[ \log \widetilde A_1\big]=  \e \big[ 1_{(\nu=1)}  A_1^{-\gamma} \, \log A_1\big]  .$

       Define  $f(t):=  \e \big [ 1_{(\nu=1)} \,   \sum_{ \vert u \vert =1} \ee^{t\, V(u)}   \big] \equiv  \e \big [ 1_{(\nu=1)} \,   A_1^{- t}\big]$ which is finite for $t \in [- \chi, \gamma]$, in particular $f(- \chi )<1$ and $f(0) <1= f(\gamma)$.       By the assumption of integrability in Theorem \ref{T:smalldev}, $\e  \big [ 1_{(\nu=1)} \,   A_1^{- \gamma } ( -\log A_1)^+\big] < \infty$ which implies that $f'(\gamma-)$ exists and equals $- \e \big[ 1_{(\nu=1)}  A_1^{-\gamma} \, \log A_1\big]$.  By convexity, $f'(\gamma-) \ge {f(\gamma)- f(0) \over \gamma} >0$. Hence \begin{equation}\label{mean1}  \e \Big[ \log \widetilde A_1\Big]= - f'(\gamma-)<0.\end{equation}

     Let $(\widetilde  A_i)_{i\ge 2}$ be a sequence of i.i.d. copies of $\widetilde A_1$ and define  $X_j:= - \sum_{i=1}^j \log \widetilde A_i $ for all $j \ge 1$.  Let $r>1$ and put  \begin{equation} \label{alphar}  \alpha_r:= \inf\{ j\ge1:  X_j  > \log r\},  \end{equation}
     
     \noindent    which is a.s. finite thanks to \eqref{mean1}.  Going back to \eqref{ite1*}, we get that \begin{eqnarray*} g(r) &\ge& \e \Big[ g(r \widetilde A_1) 1_{( r  \widetilde A_1 <   1)} \Big] +\e \Big[ g(r \widetilde A_1) 1_{( r  \widetilde A_1 \ge   1)} \Big]  \\
     &\ge&  \e \Big[ g(r \widetilde A_1) 1_{( r  \widetilde A_1 <  1)} \Big] +\e \Big[ g(r \widetilde A_1 \widetilde  A_2) 1_{( r  \widetilde A_1\ge   1)} \Big]    ,
     \end{eqnarray*}

    \noindent where to get  the last inequality,  we have applied    \eqref{ite1*}  with  $t$ replaced by $r \widetilde A_1$  and  $\widetilde A_1$ replaced by $\widetilde A_2$.   Then we obtain that    \begin{eqnarray*} g(r)  
     &\ge&  \e \Big[ g(r \widetilde A_1) 1_{( r  \widetilde A_1 <   1)} \Big] +\e \Big[ g(r \widetilde A_1 \widetilde  A_2) 1_{( r  \widetilde A_1 \ge   1,  r  \widetilde A_1 \widetilde A_2 < 1 )} \Big]   +    \e \Big[ g(r \widetilde A_1 \widetilde  A_2) 1_{( r  \widetilde A_1 \ge   1,  r  \widetilde A_1 \widetilde A_2  \ge 1 )} \Big]   \\
     &=&       \e \Big[ g(r  \prod_{i=1}^{\alpha_r} \widetilde A_i) 1_{( \alpha_r \le 2 )} \Big] +   \e \Big[ g(r \widetilde A_1 \widetilde  A_2) 1_{(\alpha_r >2 )} \Big]  . \end{eqnarray*}

By induction, we get that for any $n\ge1$, $$ g(r) \ge  \e \Big[ g(r  \prod_{i=1}^{\alpha_r} \widetilde A_i) 1_{( \alpha_r \le n )} \Big] +   \e \Big[ g(r \prod_{i=1}^n\widetilde A_i) 1_{(\alpha_r >n )} \Big] \ge \e \Big[ g(r  \prod_{i=1}^{\alpha_r} \widetilde A_i) 1_{( \alpha_r \le n )} \Big].$$ 

\noindent Since $\alpha_r <\infty$ a.s., we let $n \to \infty$ and deduce from the monotone convergence theorem  that $$  g(r) \ge  \e \Big[ g(r  \prod_{i=1}^{\alpha_r} \widetilde A_i)  \Big] = \e \Big[ g( \ee^{ - {\cal R}_r})\Big], $$

\noindent where ${\cal R}_r:= X_{\alpha_r}- \log r >0$ denotes the overshoot of the random walk $(X_j)$ at the level $\log r$.  Note  that for any $0< t\le 1$, $g(t) = t^\gamma \Phi(t) \ge \Phi(1) t^\gamma$, hence \begin{equation} \label{low1}  g(r) \ge \Phi(1) \, \e \Big[  \ee^{ -\gamma  {\cal R}_r} \Big] , \qquad \forall \, r >1. \end{equation}

     By the assumption \eqref{hyp-exp}, $\e \big[ (( - \log \widetilde A_1)^+)^2\big] = \e \big [ 1_{(\nu=1)}\,      \sum_{ \vert u \vert =1}  (V(u)^+)^2\,  \ee^{\gamma\, V(u)}   \big]  <\infty $, then by Lorden \cite{Lorden}, Theorem 1, $\sup_{r \ge1} \e \big[ {\cal R}_r\big] < \infty$.  Consequently   for some positive constant $C$, $$ g(r) \ge \Phi(1)\, \ee^{- \gamma\, \e [{\cal R}_r]}\ge C>0, \qquad \forall \, r >1.$$  
     
    \noindent Hence \begin{equation} \label{lowerPhi(r)} \Phi(r) \ge C \, r^{- \gamma}, \qquad \forall \, r >1,\end{equation}

    \noindent which implies the lower bound in \eqref{PhiPhi}.

    To prove  the upper bound in \eqref{PhiPhi},  let $a>\gamma$ be as in \eqref{hyp-exp} such that $\e [ \sum_{i=1}^\nu A_i^{- a} ] \equiv \e[ \sum_{|u|=1} \ee^{a V(u)}]<\infty.$ Choose (and then fix)  $0< \varepsilon < {1\over2} \min( a- \gamma , \gamma) $   small and   $b:= {\gamma+\varepsilon\over 2} < \gamma $.    By Lemma \ref{L:liu}, $\Phi(t) \le c_4 \,  t^{-  b}$ for all $t \ge 1$ (with $c_4 \ge1$). Since $  \Phi(t) \le 1$ for all $0<t<1$, we obtain immediately  that  \begin{equation}\label{upgt} g(t) \le  c_4\,  t^{ \gamma-b}, \qquad \forall t>0.\end{equation}

    By \eqref{Phi} and using again the notation $\widetilde A_i, i\ge 1$, we get  that for any $t>0$, \begin{eqnarray} g(t) &\le& t^\gamma \e \Big[ \Phi(t A_1) 1_{(\nu=1)}\Big] + t^\gamma \, \e \Big[ 1_{( \nu \ge2)} \Phi(t A_1) \Phi(t A_2)\Big]  \nonumber \\
    &= &  \e  \Big[  g(t \widetilde A_1)  \Big]  +  t^{-\gamma} \e \Big[ 1_{( \nu \ge2)} g(t A_1) g(t A_2) A_1^{-\gamma} A_2^{-\gamma}\Big] \nonumber \\
	&\le&  \e  \Big[  g(t \widetilde A_1)  \Big]  + c^2_4\, t^{ \gamma - 2 b}\, \e \Big[ 1_{(\nu\ge 2)} A_1^{-b} A_2^{-b} \Big]  \qquad \mbox{(by \eqref{upgt})}
	\nonumber \\&=:&  \e  \Big[  g(t \widetilde A_1)  \Big] + C_\varepsilon\, t^{-\varepsilon} , \label{gt>}
    \end{eqnarray}
    
    \noindent with $C_\varepsilon:= c^2_4\, \e \Big[ 1_{(\nu\ge 2)} A_1^{-b} A_2^{-b} \Big]  \le c^2_4  \e \Big[ \sum_{i=1}^\nu A_i^{-2b}\Big]$ by Cauchy-Schwarz' inequality. Then $C_\varepsilon<\infty$ by     the assumption \eqref{hyp-exp} and  the choice that  $b < a/2$.  
    
    Let $r>1$. As before, we shall iterate \eqref{gt>} up to the stopping time $\alpha_r$ (cf. \eqref{alphar}).  We have that \begin{eqnarray*} g(r) &\le&  C_\varepsilon r^{-\varepsilon} + \e \Big[ g( r \widetilde A_1) 1_{(\alpha_r=1)}\Big] +   \e \Big[  1_{(\alpha_r>1)} \big( C_\varepsilon (r \widetilde A_1)^{-\varepsilon} + g( r \widetilde A_1 \widetilde A_2) \big) \Big]   \\
    &=& C_\varepsilon r^{-\varepsilon} + C_\varepsilon \e \big[ (r \widetilde A_1)^{-\varepsilon} 1_{(\alpha_r>1)} \big] +  \e \Big[ g( r \prod_{i=1}^{2 \wedge \alpha_r} \widetilde A_i)  \Big]  . 
    \end{eqnarray*}    
    
    \noindent By induction, we get that for any $n \ge 2$,  \begin{eqnarray}\label{upgr1} g(r) & \le & C_\varepsilon r^{-\varepsilon} + C_\varepsilon \, \sum_{k=1}^{n-1} \e \Big[ 1_{( \alpha_r >k)} ( r \prod_{i=1}^k \widetilde A_i)^{-\varepsilon}\Big] + \e \Big[ g( r \prod_{i=1}^{n \wedge \alpha_r} \widetilde A_i)  \Big]  \nonumber \\
    &=& C_\varepsilon r^{-\varepsilon} + C_\varepsilon \e \Big[ \sum_{k =1}^{n \wedge \alpha_r-1} \ee^{ \varepsilon  (X_k- \log r)}\Big] + \e \Big[ g( r \ee^{- X_{n \wedge \alpha_r}}) \Big], \end{eqnarray}

   \noindent by using the random walk $X_j \equiv - \sum_{i=1}^j \log \widetilde A_i$, $j\ge1$. The random walk $(X_j)$ has  
    positive drift and $ \e[ X_1^2]= \e \big [ 1_{(\nu=1)}\,      \sum_{ \vert u \vert =1}  (V(u) )^2\,  \ee^{\gamma\, V(u)}   \big]  <\infty $ by the assumption \eqref{hyp-exp},  then by Lemma 5 in \cite{EOY},   $$ \e \Big[ \sum_{ k=1}^{\alpha_r - 1} \ee^{ \varepsilon  (X_k- \log r)}\Big] \le C'_\varepsilon<\infty,$$

    \noindent  for some constant $C'_\varepsilon$  independent of $r$.  On the other hand,   $g(r\ee^{-X_{\alpha_r}}) \le   1$ (since $r\ee^{-X_{\alpha_r}}\le 1$),  then we obtain that for all $r >1$, $n \ge 2$, \begin{equation}\label{upgr2} g(r) \le C_\varepsilon+C'_\varepsilon+1 + \e \Big[ g( r \ee^{- X_n}) 1_{( n < \alpha_r)}\Big]   \le C^{''}_\varepsilon + c_4\, r^ \varepsilon\, \e \Big[ \ee^{- \varepsilon X_n} 1_{( n < \alpha_r)}\Big],  \end{equation} where  in the last inequality we have  used the facts that  $t:= r \ee^{-X_n} \ge 1$ on $\{n < \alpha_r\}$ and that  $g(t)\le c_4 t^\varepsilon$  for any $t\ge1$ by   Lemma \ref{L:liu}.

    Remark that $ \e \big[ \ee^{- \varepsilon X_1}\big]=  \e \Big[  (\widetilde A_1)^\varepsilon\Big]= \e \Big[ 1_{(\nu=1)} (A_1)^{-\gamma+\varepsilon}\Big]<1$ by convexity.  Then $\e[\ee^{-\varepsilon X_n}] \to 0$ as $n \to \infty$, which in view of \eqref{upgr2} yield that for any $r>1$ ($\varepsilon$ being fixed), $g(r) \le C^{''}_\varepsilon$, i.e. $$ \Phi(r) \le C^{''}_\varepsilon \, r^{- \gamma}, \qquad \forall r>1.$$

  This and \eqref{lowerPhi(r)}  imply  \eqref{PhiPhi}: $  \Phi(r) \asymp r^{-\gamma}$ for all $r \ge1$.  The small deviation in \eqref{smalldev1new} follows from a standard Tauberian argument (see e.g. \cite{Liu99}, Lemma 4.4).      $\Box$

    \subsection{Proof of Theorem \ref{T:smalldev}:  the B\"{o}ttcher case} 
    
    The proof of \eqref{smalldev2} goes in the same spirit as that of \eqref{PhiPhi}.  Let $  h(t):= - \log \e \big[ \ee^{-t  Z } \big],  t \ge0.$   Note  that  $h$ is an increasing, concave function and vanishing at zero. Using the notations introduced in \eqref{eq-law}, we get that    $$ \ee^{- h(t)}= \e \Big[ \ee^{ - \sum_{i=1}^\nu h(  t A_i)} \Big], \qquad \forall t \ge 0.$$

 On an enlarged probability space, we may find a random variable $\xi$ such that $$ \p \Big( \xi=i \, \big \vert \, {\cal A} \Big)= { A_i^\beta \over \sum_{j=1}^\nu A_j ^\beta}, \qquad 1 \le i \le \nu,$$
 
 \noindent where ${\cal A}:= \sigma \{ A_i, 1 \le i \le \nu, \, \nu\}$.  Then  $ \sum_{i=1}^\nu h(t A_i)=  ( \sum_{i=1}^\nu A_i^\beta)\, \e \big[ { h(t A_\xi) \over A_\xi^\beta } \big\vert \, {\cal A}\big] , $  and  by Jensen's inequality, we have that  for any $t\ge0$, $$ \ee^{-  \sum_{i=1}^\nu h(t A_i)}  \le \e\left[ \exp\Big( -  ( \sum_{i=1}^\nu A_i^\beta)\,{ h(t A_\xi) \over A_\xi^\beta } \Big) \, \Big\vert\, { \cal A} \right].$$

 \noindent Write for brevity  $$ B:=  A_\xi, \qquad \eta:= {1\over A_\xi^\beta } ( \sum_{i=1}^\nu A_i^\beta)>1, \quad \mbox{a.s.}$$  
 
 \noindent  [$\eta>1$ because $\nu\ge 2$ a.s.] Then for any $t \ge0$, we have \begin{equation}\label{ite1} \ee^{- h(t)} \le \e\Big[ \ee^{ - \eta\, h(t B)}\Big].  \end{equation}

 We shall iterate  the inequality \eqref{ite1} up to some random times: Let $(\eta_i, B_i)_{i\ge1}$ be an i.i.d. copies of $(\eta, B)$. Let $r>1$ and define $$    \Upsilon_r := \inf\{i \ge 1:  \prod_{j=1}^i B_j \le {1 \over r}\}.$$
   Observe  that $$ \e \big[\log B\big]= \e \Big[ { \sum_{i=1}^\nu A_i^\beta \log A_i \over \sum_{i=1}^\nu A_i^\beta} \Big] =  - \e \Big[ { \sum_{\vert u \vert=1} \ee^{- \beta V(u)}  V(u)  \over \sum_{\vert u \vert=1} \ee^{- \beta V(u)}  } \Big]  = \psi'(\beta),$$
 
 \noindent where $\psi(x):= \e\big[ \log  \sum_{\vert u \vert=1} \ee^{- x V(u)} \big]$ for $0 \le x \le \chi$. Note that $\psi$ is convex on $[0, \chi]$, $\psi(\chi) < \log   \e\big[   \sum_{\vert u \vert=1} \ee^{-  \chi V(u)} \big] \le 0$,  and $\psi(\beta) \ge 0$ since $ \sum_{\vert u \vert=1} \ee^{- \beta V(u)} \ge 1$ by the definition of $\beta$.  By convexity, $\psi'(\beta) \le { \psi(\chi)- \psi(\beta) \over \chi- \beta} <0$. Then $\e \big[ \log B\big] < 0$ which implies that  $\Upsilon_r < \infty$, a.s.  By \eqref{ite1}, we see that for  \begin{eqnarray*} \ee^{- h(r)} &\le&  \e\Big[ \ee^{ - \eta_1\, h(r B_1)} 1_{( r B_1 \le 1)} \Big]  + \e\Big[ \ee^{ - \eta_1\, h(r B_1)} 1_{( r B_1>1)} \Big]  \\
 	&=& \e\Big[ \ee^{ - \eta_1\, h(r B_1)} 1_{( \Upsilon_r = 1)} \Big]  +  \e\Big[ \ee^{ - \eta_1\, h(r B_1)} 1_{( r B_1>1)} \Big]. 
 \end{eqnarray*}

 Applying \eqref{ite1} to $t= rB_1$, we get that $$   \ee^{ - \eta_1\, h(r B_1)}   \le  \Big( \e \big[ \ee^{- \eta_2 h( r B_1 B_2)} \, \big\vert \sigma\{\eta_1, B_1\}\big]  \Big)^{\eta_1}    \le  \e \big[ \ee^{- \eta_1 \eta_2 h( r B_1 B_2)} \, \big\vert \sigma\{\eta_1, B_1\}\big] , $$

\noindent by Jensen's inequality, since $\eta_1>1$.  It follows that 
  $  \e\big[ \ee^{ - \eta_1\, h(r B_1)} 1_{( r B_1>1)} \big] \le \e\big[  1_{( r B_1>1)}  \ee^{- \eta_1 \eta_2 h( r B_1 B_2)}  \big]$,     hence \begin{eqnarray*}  \ee^{- h(r) } &\le &\e\Big[ \ee^{ - \eta_1\, h(r B_1)} 1_{( \Upsilon_r = 1)} \Big]  + \e\Big[  1_{( r B_1>1)}  \ee^{- \eta_1 \eta_2 h( r B_1 B_2)}  \Big] \\
  	&=& \e\Big[ \ee^{ - \eta_1 \, h(r B_1)} 1_{( \Upsilon_r = 1)} \Big]  +\e\Big[ \ee^{ - \eta_1 \eta_2 \, h(r B_1B_2)} 1_{( \Upsilon_r = 2)} \Big] + \e\Big[  1_{( r B_1B_2>1)}  \ee^{- \eta_1 \eta_2 h( r B_1 B_2)}  \Big]. \end{eqnarray*}
  
  Again applying \eqref{ite1} to $t=r B_1B_2$ and using Jensen's inequality (since $\eta_1\eta_2 >1$), we get that $\e\big[  1_{( r B_1B_2>1)}  \ee^{- \eta_1 \eta_2 h( r B_1 B_2)}  \big]  \le \e\big[  1_{( r B_1B_2>1)}  \ee^{- \eta_1 \eta_2  \eta_3 h( r B_1 B_2 B_3)}  \big]$, and so on. We get that for any $n \ge 1$, \begin{eqnarray} \ee^{- h(r)} &\le& \e \Big[ \ee^{ - (\prod_{i=1}^{\Upsilon_r}  \eta_i) \, h( r \prod_{i=1}^{\Upsilon_r}  B_i) } 1_{( \Upsilon_r  \le n)} \Big] + \e \Big[ \ee^{- (\prod_{i=1}^n \eta_i )\, h( r \prod_{i=1}^n B_i) } 1_{( \Upsilon_r  > n)} \Big]  \nonumber  \\
   &=:& A_{\eqref{ite2}}+ C_{\eqref{ite2}}. \label{ite2}
  \end{eqnarray}
  
 By \eqref{bounded},  $ B \ge \ee^{-K}$ a.s., then ${ 1 \over r} \ge \prod_{i=1}^{\Upsilon_r}  B_i  > {1\over r} \ee^{- K}$. Notice that by \eqref{defbeta} the definition of $\beta$,  $\sum_{i=1}^\nu A_i^\beta \ge 1$ a.s.; Then $\eta \ge B^{- \beta}$ and  $\prod_{i=1}^{\Upsilon_r}  \eta_i \ge   r^\beta$. It follows that for any $n$,  $$ A_{\eqref{ite2}}\le \ee^{ - r^\beta\, h(\ee^{-K})}.$$

 To deal with $C_{\eqref{ite2}}$, we remark that   on $\{\Upsilon_r >n\}$, $r \prod_{i=1}^n B_i\ge 1$.  It follows that $$ C_{\eqref{ite2}} \le \e \Big[ \ee^{- h(1) \prod_{i=1}^n \eta_i  }\Big].$$
 
 \noindent Since $\eta_i >1$ a.s., $ \prod_{i=1}^n \eta_i  \uparrow \infty$ as $n \to \infty$, then by the monotone convergence theorem $ \limsup_{n \to \infty} C_{\eqref{ite2}}=0$. Letting  $n\to \infty$ in \eqref{ite2},   we   obtain that  \begin{equation}\label{uplemmabo} \e\Big[ \ee^{-r\, Z}\Big] \equiv \ee^{- h(r)} \le \ee^{- h(\ee^{-K})\, r^\beta}, \qquad \forall \,r >1,
 \end{equation} 
 
 \noindent which is stronger than the upper bound in \eqref{smalldev2}.

 To prove the lower bound,     recalling   that $\mbox{essinf } \sum_{i=1}^\nu A_i^\beta=1$   and    $A_i \ge \ee^{-K}$,  we deduce that    for any small $\varepsilon>0$,  there are   some integer $2 \le k \le  \mbox{esssup}\, \nu $,  and some real numbers $a_1, ..., a_k \in (0, 1)$ such that $  \sum_{i=1}^k  a_i ^\beta \ge 1$ and $  \sum_{i=1}^k  a_i ^{\beta +\varepsilon} < 1$ and $ p:= \p\big( A_i  \le a_i, \forall 1\le i \le k, \, \nu=k\big) >0$.  Therefore $$ \ee^{- h(t)}= \e \Big[ \ee^{- \sum_{i=1}^\nu h( t A_i)}\Big] \ge p\,  \ee^{- \sum_{i=1}^k h(t a_i)}, \qquad t \ge 0.$$ 
 
\noindent Let $b:= \log (1/p) >0$ and define a random variable $Y \in \{ a_1, ..., a_k\}$ such that for any measurable and nonnegative function $f$, $\e\big[ f(Y)\big]= { 1 \over k} \sum_{i=1}^k f( a_i)$.   Therefore,   \begin{equation}\label{ite4}  h(t) \le b + k \, \e\big[ h( t Y) \big], \qquad \forall\, t \ge 0.
 \end{equation}

As in the proof of the upper bound, we shall iterate the above inequality up to some random times: Let $(Y_j)_{j\ge1}$ be an i.i.d. copies of $Y$. For $r>1$, we define $$ \theta:=\theta_r:= \inf\{ j\ge 1: \prod_{i=1}^j Y_i \le {1\over r}\}.$$

 Since $Y \le \max_{1\le i \le k} a_i <1$, $\theta$ is a bounded random variable. Going back to \eqref{ite4}, we get that \begin{eqnarray*} h(r)  &\le& b +k \, \e \Big[ h( r Y_1) 1_{( r Y_1 \le 1)}\Big] + k\, \e \Big[ h( r Y_1) 1_{( r Y_1 > 1)}\Big] \\
	&\le& b +k\,  \e \Big[ h( r Y_1) 1_{( \theta=1)}\Big]  + k\, \e \Big[  1_{( r Y_1 > 1)}  ( b + k \, h( r Y_1 Y_2)) \Big]  \\
	&=& b + k\, \e \Big[ h( r Y_1) 1_{( \theta=1)}\Big]  +  bk\,  \p\big( \theta >1\big) +  k^2 \e \Big[  1_{( r Y_1 > 1)}   h( r Y_1 Y_2)) \Big] .
\end{eqnarray*}

\noindent By induction, we get that for any $n \ge 1$, \begin{eqnarray} h(r) & \le& b \, \sum_{j=0}^n \, k^j \, \p\big( \theta > j\big) + \e \Big[ k^{ \theta\wedge n}\, h( r \prod_{i=1}^{\theta\wedge n} Y_i)\Big] \nonumber \\
	&=:& A_{\eqref{ite5}} + C_{\eqref{ite5}} . \label{ite5}
\end{eqnarray}

\noindent Elementary computations yield that $$ A_{\eqref{ite5}}= { b \over k-1} \e \Big[ k^{ \theta\wedge (n+1)} - 1\Big]  \le { b \over k-1} \e \Big[ k^{ \theta }  \Big]    .  $$

\noindent  Recalling $\theta$ is bounded hence $ \e \Big[ k^{ \theta }  \Big]  < \infty$. For $C_{\eqref{ite5}}$, we use the fact that $Y_i \le \max_{1\le j \le k} a_j=: \overline a < 1$.  Remark that $r \prod_{i=1}^{n } Y_i \le 1$. Then  \begin{eqnarray*} C_{\eqref{ite5}} &:=&  \e \Big[ k^{ \theta }\, h( r \prod_{i=1}^{\theta } Y_i) 1_{( \theta\le n)} \Big]   +  \e \Big[ k^{ n }\, h( r \prod_{i=1}^{n } Y_i) 1_{( \theta> n)} \Big]   \\
	&\le& h(1) \e \big[ k^\theta\big] + h(r   \overline a  ^n) \, \e \big[ k^n 1_{( \theta>n)}\big] \\
	&\le& h(1) \e \big[ k^\theta\big] + h(r   \overline a ^n) \, \e \big[ k^\theta\big].
\end{eqnarray*} 

\noindent Since $r   \overline a ^n \to 0$ as $n \to \infty$, we get that [recalling that $\theta$ depends on $r$] \begin{equation}\label{uph1}  h(r) \le ( h(1) + { b \over k-1}) \, \e \big[ k^\theta\big], \qquad \forall\, r>1.
\end{equation}

To estimate $\e \big[ k^\theta\big]$, let us find $\lambda>0$ such that $ \e \Big[ Y^{\lambda}\Big]= {1\over k}.$ By the law of $Y$, this is equivalent to $  \sum_{ i=1}^k  a_i^\lambda= 1.$  

By the choice of $(a_i)$, we have $\beta \le \lambda < \beta+\varepsilon$.  Then the process $n \to k^n\, \prod_{i=1}^n Y_i^\lambda $ is a martingale (moreover uniformly integrable on $[0, \theta]$). Hence the optional stopping   theorem implies that $$ 1= \e \Big[ k^\theta \prod_{i=1}^\theta Y_i^\lambda \Big] \ge \e \Big[ k^\theta  \Big] \,  r^{-\lambda}\, \min_{ 1 \le i \le k } a_i^\lambda,$$ 

\noindent since $\prod_{i=1}^\theta Y_i \ge {1\over r} \min_{ 1 \le i \le k } a_i$.  This and \eqref{uph1} give that $$ h(r) \le  ( h(1) + { b \over k-1}) \max_{ 1\le i \le k} a_i^{-\lambda}\, r^\lambda, \qquad \forall r >1,$$

\noindent yielding the lower bound in \eqref{smalldev2}  since $\lambda < \beta+\varepsilon$. This  completes  the proof of \eqref{smalldev2}.  Finally,  by using the elementary inequalities:  for any $\varepsilon, t >0$, $\ee^{- \varepsilon t } \p( Z < \varepsilon) \le \e [ \ee^{- t Z}] \le \p (Z < \varepsilon) + \ee^{- \varepsilon t}, $   we  immediately      deduce from \eqref{smalldev2} that $\p ( Z < \varepsilon) = \ee^{- \varepsilon^{-\beta/(1-\beta) +o(1)}}$ as $\varepsilon \to 0$.  $\Box$.

 \section{Proof of Theorem \ref{t:3log}} \label{S:4}

Let us  give some preliminary estimates on the branching random walk:

   \begin{lemma}\label{L:infsup} Assume \eqref{hyp} and \eqref{int1}. There exists some constants $c_5, c_6>0$ such that for $n \ge 1$,   \begin{equation}\label{infsup} \p\Big( \min_{ \vert u \vert =n} \overline V(u) <   c_5 \,   n^{1/3} \Big) \le c_6\, \ee^{ -c_5 n^{1/3}},\end{equation} where we recall that  for any $\vert u \vert =n$, $\overline V(u):= \max_{ 1\le i \le n} V(u_i)$.  Consequently, for any $ 0<\lambda \le c_5 n^{1/3}$, we have \begin{equation} \label{infsup2} \p \Big( \max_{u \in \L_\lambda} \vert u \vert > n \Big) \le  c_6 \, \ee^{- c_5\, n^{1/3}}. \end{equation} 
  \end{lemma}
  
 We mention that under an  extra integrability condition, i.e. $\exists \delta >0$ such that $\e [\nu^{1+\delta}] < \infty$,  $ n^{-1/3} \min_{ \vert u \vert =n} \overline V(u) \to ({ 3 \pi^2 \sigma^2\over 2})^{1/3}$  $\p^*$-a.s. (see \cite{FHS} and \cite{FZ10}) and the probability term in \eqref{infsup} is equal to $\ee^{ (c_5 - ({ 3 \pi^2 \sigma^2\over 2})^{1/3} +o(1)) n^{1/3}}$ for any $0 < c_5 < ({ 3 \pi^2 \sigma^2\over 2})^{1/3}$ (see \cite{FHS}, Proposition 2.3).   Here, we  only assume   \eqref{hyp} and \eqref{int1}, and we do not seek the precise upper bound in \eqref{infsup}.

 \medskip
 {\noindent \bf Proof of Lemma \ref{L:infsup}.}  We shall use the following fact (see Shi \cite{Shi}): \begin{equation}\label{estimate-infV}  \p( \inf_{ u \in \T} V(u) < - \lambda\Big) \le \ee^{- \lambda}, \qquad \forall \, \lambda \ge0. \end{equation} 
  
  \noindent 
   Consider $0< c< ({ \pi^2 \sigma^2 \over 8  })^{1/3} $.  Then \begin{eqnarray*}    \p\Big( \min_{ \vert u \vert =n} \overline V(u) < c n^{1/3} , \inf_{ u \in \T} V(u) \ge  -  c n^{1/3}  \Big)  &\le & \e \Big[ \sum_{ \vert u \vert =n} 1_{(\max_{ 1 \le i \le n} \vert V(u_i) \vert \le c n^{1/3})}\Big] \\
   &=& \e \Big[ \ee^{ S_n} 1_{( \max_{ 1 \le i \le n} \vert S_i\vert \le c n^{1/3})} \Big]  \qquad \mbox{(by \eqref{many})}\\
   &\le & \ee^{ c n^{1/3}} \, \p \Big(  \max_{ 1 \le i \le n} \vert S_i\vert \le c n^{1/3}\Big)  \\
   &=& \ee^{ c n^{1/3}} \, \ee^{- ({ \pi^2 \sigma^2 \over 8 c^2}+o(1) ) n^{1/3}},
   \end{eqnarray*}
   
   \noindent where the last equality follows from    Mogulskii \cite{M74}.   This and    \eqref{estimate-infV}  easily yield      the Lemma by choosing a sufficiently small constant $c$. 
$\Box$

\medskip

  Recall \eqref{defll}. Define for $a \in (0, \infty]$ and $\lambda>0$,  \begin{equation} \label{La} \L_\lambda^{(a)}:= \Big\{ u \in \L_\lambda: V(u)  \le  \lambda+ a \Big\}. \end{equation}

\noindent In particular,   $\L_\lambda^{(\infty)}=\L_\lambda$.      Recall \eqref{DL}.  Since the function $x \to x \ee^{-x}$ is decreasing for $x \ge 1$, then for any $\lambda >1$,  $D_{\L_\lambda} \le   \lambda \ee^{- \lambda} \# \L_\lambda$, which implies that    \begin{equation} \label{lowb}  \liminf_{ \lambda \to \infty} \lambda \ee^{- \lambda}\, \# \L_\lambda \ge D_\infty >0.  \qquad \mbox{a.s. on $\S$.} \end{equation}

If $\nu=\infty$ [which is allowed under \eqref{hyp} and \eqref{int1}], then $\#\L_\lambda=\infty$ hence \eqref{lowb} cannot be strengthened into a true limit.   We present  a similar result   for $\L_\lambda^{(a)}$:

\begin{lemma} \label{L:Nerman} Assume \eqref{hyp}, \eqref{int1} and that  $\e \big[ \sum_{\vert u \vert=1} (V(u)^+)^3 \ee^{- V(u)}\big] < \infty$.  There exists some $a_0>0$ such that  for all large  $a\ge a_0$, almost surely on the set of non-extinction $\S$, $$ 0<  \liminf_{\lambda \to \infty} \lambda \,\ee^{-\lambda} \,  \#\L_\lambda^{(a)} \le \limsup_{\lambda \to \infty} \lambda \,\ee^{-\lambda} \,  \#\L_\lambda^{(a)} < \infty.  $$\end{lemma}

{\noindent\bf Proof of Lemma \ref{L:Nerman}.}   We only deal with the case when   the distribution  of $\Theta$ is non-lattice, in this case,  the limit exists. The lattice case can be treated  in a similar way, by applying  Gatzouras (\cite{G00}, Theorem 5.2),  a lattice version of Nerman \cite{N81}'s  result, but the cyclic phenomenon could prevent from the existence of limit.  In the non-lattice case, we are going to prove that   for any $a>0$, almost surely on the set of non-extinction $\S$,  \begin{equation}\label{Nerman2}   \lim_{\lambda \to \infty} \lambda \ee^{-\lambda} \,  \#\L_\lambda^{(a)}=  c_7(a) \, D_\infty, \end{equation}
where $c_7(a):=  {  1 \over \e\big[ S_{\tau_0}\big]}  \e \big[ \ee^{ \min(a,S_{\tau_0})}- 1\big]  $, and $S_\cdot$  and $\tau_0$   are defined by      \eqref{many} and  \eqref{tau0} respectively.  Obviously, $c_7(a)>0$ for all large $a$.

To get  \eqref{Nerman2},  we consider a new  point process $\widehat  \Theta:=\sum_{u \in \L_0} \delta_{\{ V(u)\}}$ on $( 0, \infty)$.   Generate a branching random walk $(\widehat V(u),  u \in \widehat \T)$ from the point process $\widehat \Theta$, in the same way as $(V(u), u \in \T)$ do from $\Theta$. 
 Remark that  $\S= \{\sup_{u \in \T} V(u) =\infty\}=\{ \widehat \T \mbox { is infinite}\}$,   and $$ \#\L_\lambda^{(a)}= \sum_{ u \in \widehat \T} \phi_u( \lambda- \widehat V(u)), \quad \sum_{u \in   \L_\lambda} \ee^{- V(u) +\lambda}= \sum_{ u \in \widehat \T} \psi_u(\lambda- \widehat V(u)), $$

 \noindent where $$ \phi_u(y):= 1_{(y\ge0)} \sum_{ v: {\buildrel \leftarrow \over v} =u} 1_{( y  < \widehat V(v)- \widehat V(u) \le  y+ a)}, \quad \psi_u(y):= 1_{(y\ge0)} \sum_{ v: {\buildrel \leftarrow \over v} =u} \ee^{y- ( \widehat V(v)- \widehat V(u))} 1_{(  \widehat V(v)- \widehat V(u)  > y)}.$$

 \noindent Applying Theorem 6.3 in Nerman \cite{N81} (with $\alpha=1$   there) gives that almost surely on  $  \S$, $$ {  \sum_{ u \in \widehat \T} \phi_u( \lambda- \widehat V(u)) \over   \sum_{ u \in \widehat \T} \psi_u( \lambda- \widehat V(u))} \, \to \, { \e \big[ \sum_{ \vert u \vert=1, u \in \widehat \T} ( \ee^{- (\widehat V(u) - a)^+}- \ee^{- \widehat V(u)}) \big]\over \e \big[ \sum_{ \vert u \vert=1, u \in \widehat \T}  \widehat V(u) \ee^{- \widehat V(u)}\big]}.$$

 Remark that $\e \big[ \sum_{ \vert u \vert=1, u \in \widehat \T} ( \ee^{- (\widehat V(u) -a)^+}- \ee^{- \widehat V(u)}) \big]= \e \big[ \sum_{u \in \L_0} (\ee^{- (V(u) - a)^+}- \ee^{-V(u)})\big]= \e \big[ \ee^{ \min(a,S_{\tau_0})}- 1\big]$ and $\e \big[ \sum_{ \vert u \vert=1, u \in \widehat \T}  \widehat V(u) \ee^{- \widehat V(u)}\big]= \e \big[ \sum_{ u \in \L_0}    V(u) \ee^{-   V(u)}\big] = \e\big[ S_{\tau_0}\big]$.  Hence on $\S$, a.s., \begin{equation} \label{Nerman3}  { \#\L_\lambda^{(a)} \over \sum_{u \in   \L_\lambda} \ee^{- V(u) +\lambda}} \to  c_7(a).\end{equation}
 
 On the other hand,  almost surely, \begin{equation}\label{ddd4} D_{\L_\lambda}= \lambda \ee^{-\lambda} \Big( \sum_{ u \in \L_\lambda} \ee^{- V(u) +\lambda}  + {1\over \lambda} \eta_\lambda\Big) \, \to \, D_\infty, \qquad \lambda\to \infty, \end{equation}
 
 \noindent where $\eta_\lambda:= \sum_{u \in \L_\lambda} ( V(u) - \lambda) \ee^{- V(u)+ \lambda}$.  By the many-to-one formula and the assumption, $\e \big[  (S_1^+)^3\big] =\e \big[ \sum_{\vert u \vert=1} ( V(u)^+)^3 \ee^{- V(u)}\big] < \infty$. Then by Doney \cite{D80}, $\e \big[ S_{\tau_0}^2\big] < \infty$. 
 
 Note that $\eta_\lambda=   \sum_{ u \in \widehat \T} \widetilde \psi_u(\lambda- \widehat V(u))$ with $ \widetilde \psi_u(y):= 1_{(y\ge0)} \sum_{ v: {\buildrel \leftarrow \over v} =u} \ee^{y- ( \widehat V(v)- \widehat V(u))} (\widehat V(v)- \widehat V(u)-y)1_{(  \widehat V(v)- \widehat V(u)  > y)}$. 
   In the same manner  we get  that   almost surely on $\S$,  \begin{equation}  \label{es-eta}  \lim_{\lambda\to\infty}{ \eta_\lambda\over\sum_{ u \in \L_\lambda} \ee^{- V(u) +\lambda}  } =c_8,  \end{equation}

 \noindent
 with $c_8:= {1\over2} \, { \e \big[ S_{\tau_0}^2 \big] \over \e\big[ S_{\tau_0}\big]}>0$. It follows that a.s. on $\S$, $\sum_{ u \in \L_\lambda} \ee^{- V(u) +\lambda}   \sim {1\over \lambda} \,\ee^\lambda D_{\L_\lambda} \sim  {1\over \lambda} \ee^\lambda D_\infty   $ as $\lambda \to \infty$. This combined with \eqref{Nerman3} and \eqref{ddd4} yield \eqref{Nerman2}, as desired. $\Box$

 \begin{remark}\label{remark-nerman}  The condition $\e \big[ \sum_{\vert u \vert=1} (V(u)^+)^3 \ee^{- V(u)}\big] < \infty$   was       used    in the  above  proof of Lemma \ref{L:Nerman}  only to obtain \eqref{es-eta} which controls   the contribution of $\eta_\lambda$ in $D_{\L_\lambda}$.  We do  not know how to relax this condition.
 \end{remark}

  We consider now  some deviations on the minimum $\M_n$.  If the distribution of $\Theta$ is non-lattice, A{\"{\i}}d{\'e}kon (Proposition 4.1,  \cite{A11})   proved   that   for any $A>0$ and  for all large $n, \lambda$  such that  $A \le \lambda \le  {3\over2} \log n -A$,  $$ \p \Big( \M_n < {3 \over 2} \log n - \lambda \Big) = (c_9 +o_A(1)) \, \lambda \, \ee^{-\lambda},$$ with $c_9$ some positive constant and $o_A(1) \to 0$  as $A\to \infty$ uniformly in $n, \lambda$.  We shall need in the proof of  Theorem \ref{t:3log} an estimate which holds uniformly in $\lambda$.   
  
   \begin{lemma}[Mallein \cite{Mallein}, Lemma 4.2]  \label{L:d1} Assume \eqref{hyp} and \eqref{int1}. There is some constant $c_{10}>0$ such that   $$ \p \Big( \M_n <  {3 \over 2} \log n - \lambda \Big) \le c_{10}\, (1+\lambda) \ee^{-\lambda} , \qquad \forall n \ge 1, \, \lambda\ge0.$$
   \end{lemma}

         \medskip

  
  The tightness of $(\M_n - {3\over2} \log n)_{n\ge1}$ under \eqref{hyp} and \eqref{int1} was implicitly contained in A{\"{\i}}d{\'e}kon (\cite{A11}) (see also \cite{BZ06}, and see \cite{AbR} for exponential decay  under some additional assumptions): Assume \eqref{hyp} and \eqref{int1}. We have\footnote{In fact, by Lemma 3.6 in \cite{A11} and using the fact that $\M_n$ is stochastically smaller than $M^{\mathrm {kill}}_n$,  we obtain that $\sup_{n\ge3} \p(\M_n\ge {3\over 2} \log n) \le \ee^{-C}$  for  some (small) constant $C>0$. For any $k\ge1$, denote by $Z_k:=\sum_{|u|=k}1$ the number of individuals at generation $k$. By the triangular inequality and the branching property at $k$, we get that  for any $n\ge k+3$, $\p\big(  \M_n \ge  { 3\over 2} \log n   +  \lambda, \, \S\big) \le \p\big( \exists |u|=k: V(u) > \lambda\big) + \e \big[ 1_{(Z_k>0)} \ee^{-C Z_k}\big]$. Letting $\lambda\to\infty$ and then $k\to \infty$, we get   \eqref{AR-new}. The left tail $\limsup_{n\to\infty}\p\big(  \M_n -  { 3\over 2} \log n  < - \lambda \big) $, as $\lambda\to  \infty$,  follows from  \cite{A11}, see also Lemma \ref{L:d1}.}   \begin{equation}\label{AR-new}  \limsup_{\lambda\to\infty}\limsup_{n\to\infty}  \p^* \Big(   \M_n -  { 3\over 2} \log n   \ge  \lambda\Big) =0,\end{equation}

 \noindent where as before, $\p^*(\cdot):= \p( \cdot \vert \S)$.    We need some tightness uniformly in $n$: 
 
\begin{lemma}\label{tightness} Assume \eqref{hyp} and \eqref{int1}. For any fixed  $a>1$, we have $$ \limsup_{ n \to \infty} \, \p^* \Big( \max_{ n \le k \le a n} \M_k \ge {3\over2} \log n + x \Big) \to 0, \qquad \mbox{ as } x \to \infty.$$ 
\end{lemma}

 \medskip
  {\noindent \bf Proof of Lemma \ref{tightness}.} Obviously, it is enough to prove the Lemma for $a=2$.   By Lemma \ref{L:d1}, there exists some $\lambda_0>0$ such that for all $\lambda\ge \lambda_0$ and for all $n \le k \le 3n$, \begin{equation}\label{k3nlambda} \p \Big( \M_{ 4n -k} \ge {3\over2}\log n - \lambda\Big) \ge \exp\Big( - 2 c_{10}  \, \lambda \, \ee^{- \lambda} \Big). \end{equation}

  Let $x\ge 2 \lambda_0$ and $n \gg x$.  Define  $$ \kappa_x\equiv \kappa_x(n):= \inf\{ k \ge n: \M_k \ge {3 \over 2}\log n + x \}, \qquad (\inf_\emptyset =\infty).$$ 
 
Let $n \le k \le 3n$.   Denote by $\S_k$  the event that  the Galton-Watson tree $\T$ survives up to the generation $k$.   Then $\S_k$ is non-increasing  on $k$. On the set $\{ \kappa_x= k\} \cap \S_k $, $V(u) > {3\over 2} \log n +x$ for any $\vert u \vert =k$. Let $0< y < x-\lambda_0$.   It follows from the branching property that  on $\{ \kappa_x= k\} \cap  \S_k $,  \begin{eqnarray*}   \p\Big( \M_{4n} > {3 \over2} \log n + y \,  \big\vert\, \F_k  \Big)  
 	&=&   \prod_{\vert u \vert =k} \p \Big(   \M_{4n - k} \ge {3\over2} \log n - \lambda \Big)\big\vert_{\lambda= V(u)- y } \\	&\ge& \exp \Big( - 2 c_{10}\, \sum_{ \vert u \vert =k} ( V(u) - y) \ee^{- (V(u)- y)} \Big) \\
	& \ge &  \exp \Big( - 2 c_{10}\,  \ee^y  D_k   \Big) ,  	 \end{eqnarray*}

\noindent where we have used  \eqref{k3nlambda} to get the above first  inequality.

	 Therefore for any $\varepsilon>0$ and $n \le k \le 3n$, \begin{eqnarray}    \p\Big( \M_{4n} > {3 \over2} \log n + y , \, \S_k,  \, \kappa_x= k  \Big)  		& \ge & \e \Big[   \ee^{-  2  c_{10}\,   \ee^y\, D_k} 1_{( \S_k \cap\{  \kappa_x=k\})}\Big] \nonumber
		\\ &\ge& \ee^{-\varepsilon}\, \p \Big(A_{\eqref{4ny4}}\, \cap\, \{ \kappa_x=k\}\Big), \label{4ny2}
		\end{eqnarray}
	
	\noindent where \begin{equation} \label{4ny4} A_{\eqref{4ny4}}:=\S\cap  \big\{   \sup_{ n \le j \le 4n}   D_j  \le {\varepsilon\over 2c_{10}} \ee^{-y}\big\}.\end{equation}

	\noindent  Since $\S_k \subset \S_n$ for $k \ge n$, \eqref{4ny2} still  holds if we replace $\S_k$ by $\S_n$ in the LHS.  Taking the sum over $n \le k\le 3 n $  for \eqref{4ny2} (with  $\S_k$ replaced  by $\S_n$), we get that for any $\varepsilon>0$, $0< y <x-\lambda_0$ and all $n \ge n_0$,   \begin{eqnarray}  \p\Big( A_{\eqref{4ny4}} \cap \big\{  \max_{ n \le k \le 3 n} \M_k \ge {3 \over 2} \log n + x\big\}\Big) &\le& \ee^{ \varepsilon} \,  \p\Big( \M_{4n} > {3 \over2} \log n + y , \, \S_n   \Big)  \nonumber \\
	&\le&    \ee^{ \varepsilon} \,  \p\Big( \M_{4n} > {3 \over2} \log n + y ,  \S  \Big)   +  \ee^\varepsilon  \p(\S^c \cap \S_n)	\nonumber \\
	&\le&  \varepsilon + \ee^\varepsilon  \p(\S^c \cap \S_n), \label{4ny3b}  \end{eqnarray} 
	
	\noindent by using \eqref{AR-new} if we choose a sufficiently large constant   $y= y(\varepsilon)$ only depending on $\varepsilon$.   Since  $\lim_{n \to \infty} \p(\S^c \cap \S_n)=0$,    then  for $x >  y(\varepsilon) + \lambda_0$ and all large $n \ge n_1(\varepsilon)$, \begin{equation} \label{4ny5}  \p\Big( A_{\eqref{4ny4}} \cap \big\{  \max_{ n \le k \le 3 n} \M_k \ge { 3 \over 2} \log n + x\big\} \Big) \le  2 \varepsilon. 
\end{equation}

\noindent   Note the factor $3 n$ in the above estimate and we fix our choice of $y\equiv y(\varepsilon)$ in $A_{\eqref{4ny4}}$.

Now, we shall get rid of the term $A_{\eqref{4ny4}}$ in \eqref{4ny5}.  Let $z \in (y , \, x-\lambda_0)$. Recalling the definition of $\L_z$  in \eqref{defll}.   Define \begin{equation}\label{2n2x}
A_{\eqref{2n2x}}:= \Big\{ \exists u \in \L_z: \vert u \vert \le x, V(u) \le  x ,  \, \sup_{\frac{n}{2}\le j\le 3 n }D_j^{(  u)} \le  {\varepsilon\over 2c_{10}} \ee^{-y} ,   \S^{(u)}\Big\},
\end{equation}


\noindent where   $(D_j^{(u)}, j\ge 0), \M_\cdot^{(u)}, \S^{(u)}$ are defined from the subtree $\T_u$ in the same way as $(D_j, j\ge 0), \M_\cdot, \S$ do from $\T$.   Let $n >2 x$.    The event   $\{\max_{ n \le k \le 2n} \M_k \ge  {3\over2} \log n + 2x, \S\}$ implies that for some $n \le k \le 2n$,  for any $\vert v \vert =k$, $V(v) \ge {3\over2} \log n + 2x$. If  $A_{\eqref{2n2x}} \neq \emptyset$, then we take an arbitrary  $u \in A_{\eqref{2n2x}}$ and  get that $\M^{(u)}_{k- \vert u \vert } \ge {3\over2} \log n + 2x  - V(u) \ge {3\over2} \log n +  x. $  By conditioning on $|u|$, we get  that for $ z\in (   y, \, x-\lambda_0)$ and for all   large $n\ge n_2(x, \varepsilon)$,  
 \begin{eqnarray}  &&  \p\Big( \max_{ n \le k \le 2n} \M_k \ge  {3\over2} \log n + 2x, \S , \,  A_{\eqref{2n2x}} \neq \emptyset \Big)  \nonumber\\ &\le&   	
 	\max_{1\le i \le x} \, \p\Big(    \max_{ n-i \le k \le  2n-i } \M_k \ge { 3 \over 2} \log n +  x\, , \sup_{\frac{n}{2}\le j\le 3 n }D_j \le  {\varepsilon\over 2c_{10}} \ee^{-y} ,   \S\Big)    \nonumber \\
	&\le&
	 \p\Big(    \max_{ \frac{2n}{3} \le k \le  2n } \M_k \ge { 3 \over 2} \log n +  x\, , \sup_{\frac{n}{2}\le j\le 3 n }D_j \le  {\varepsilon\over 2c_{10}} \ee^{-y} ,   \S\Big)
	 \nonumber\\
	& \le&  2 \varepsilon  , \label{2n2x2}  \end{eqnarray} 

\noindent by applying \eqref{4ny5}  to get the last inequality.    On the other hand, \begin{eqnarray}   && \p\Big( A_{\eqref{2n2x}} = \emptyset , \S\Big)   \nonumber 
	\\ &\le&  \p \Big( \forall u \in \L_z, \sup_{\frac{n}{2}\le j\le 3 n }D_j^{( u)} \ge {\varepsilon \over 2 c_{10}} \ee^{-y} \mbox{ or } (\S^{(u)})^c, \, \L_z\neq \emptyset\Big) + \p\Big( \max_{u \in \L_z } \max(V(u), \vert u \vert) \ge  x \Big)  \nonumber \\
	&=&  \e \Big[ \ee^{- p_n(\varepsilon,y)\#\L_z} 1_{( \# \L_z >0)}\Big] + \p\Big( \max_{u \in \L_z } \max(V(u), \vert u \vert) \ge  x \Big), \label{3n3x3}
\end{eqnarray}

 \noindent where the last equality is due to     the branching property at $\L_z$, and   $p_n(\varepsilon, y)$ is defined by  $\ee^{-p_n(\varepsilon, y)}:= \p\big( \sup_{\frac{n}{2}\le j\le 3 n }D_j \ge  {\varepsilon \over 2c_{10}} \ee^{-y} \mbox{ or } \S^c\big)$. As $D_n \to D_\infty$, almost surely on $\S$, we see that $\limsup_{n\to \infty} \ee^{-p_n(\varepsilon, y)}\le  \p\big( D_\infty >  {\varepsilon \over 2c_{10}} \ee^{-y} \mbox{ or } \S^c\big)$ as $n \to \infty$.  
 
 By Liu \cite{Liu01}, Theorem  2.6 \footnote{This is where we use the condition that the law of $\log \sum_{|u|=1} \ee^{-V(u)}$ is non-lattice.}, $\p\big( 0<D_\infty \le  {\varepsilon \over 2c_{10}} \ee^{-y}, \S) >0$ hence $\p\big( D_\infty >  {\varepsilon \over 2c_{10}} \ee^{-y} \mbox{ or } \S^c\big)<1$ and there exists some $p(\varepsilon,y)>0$ such that for all large $n$, $\ee^{-p_n(\varepsilon, y)}\le \ee^{- p(\varepsilon,y)}$.  

Assembling  \eqref{2n2x2} and \eqref{3n3x3}  give that  for any $z> y$,   \begin{eqnarray}  C_{\eqref{xn2n}}&:=&   \limsup_{x \to \infty} \limsup_{n \to \infty} \p\Big(   \max_{ n \le k \le 2n} \M_k \ge  {3\over2} \log n + 2x, \S\Big)   \nonumber
	\\& \le &    \e \Big[ \ee^{- p(\varepsilon,y)\#\L_z} 1_{( \# \L_z >0)}\Big] + 2 \varepsilon, \label{xn2n}  \end{eqnarray} 

Notice that $\{\#\L_z >0\}$ is nonincreasing on $z$ and its limit as $z \to \infty$ equals $\S$. Then $\p\big( \{ \#\L_z>0\} \cap \S^c\big) \to0$ as $z \to \infty$. On $\S$, we have from \eqref{lowb} that $\L_z \to \infty$  as $z \to \infty$ almost surely, hence     $\e \big[ \ee^{- p(\varepsilon,y)\#\L_z} 1_{( \# \L_z >0)}\big] \le   \e \big[ \ee^{- p(\varepsilon,y)\#\L_z} 1_{\S}\big] +  \p\big( \{ \#\L_z>0\} \cap \S^c\big) \to 0$ as $z\to\infty$.    Then letting $z\to \infty $, we see that $ C_{\eqref{xn2n}} \le  2\varepsilon$. This proves the Lemma since $\varepsilon$ can be arbitrarily small.   $\Box$

 \medskip
  
 We are now ready to give the proof of Theorem \ref{t:3log}.
 \medskip
 
{\noindent \bf Proof of   Theorem \ref{t:3log}.}  

 {\noindent \it Proof of the lower bound in Theorem \ref{t:3log}.} Consider large integer $j$.  Let $n_j:= 2^j $ and $\lambda_j:= a \log \log \log n_j$ with some constant $0<a<1$. Fix $\alpha>0$ and  put  $$A_j:= \Big\{ \M_{n_j} > {3\over2} \log n_j + \lambda_j \Big\}.$$

 \noindent   Recall  that  if the system dies out at generation $n_j$, then by definition $\M_{n_j}= \infty$. Define  $  \M_\cdot^{(u)} $     from the subtree $\T_u$ in the same way as $  \M_\cdot $ does from $\T$. Then  $A_j=\{ \forall \, \vert u \vert = n_{j-1}$, $ \M^{(u)}_{n_j - n_{j-1}} \ge  {3\over2} \log n_j + \lambda_j - V(u)\}$, which  by the branching property at $n_{j-1}$ implies that  \begin{eqnarray*}  \p\Big( A_j  \, \vert\, \F_{n_{j-1}}\Big) 	 
 	&=&     \prod_{  \vert u \vert = n_{j-1}} \p\Big(  \M_{n_j - n_{j-1}} \ge  {3\over2} \log n_j + \lambda_j - x \Big)\big\vert_{x=V(u)},  
 \end{eqnarray*}

 \noindent with the usual convention: $\prod_\emptyset:=1$.   By the lower limits of $\M_{n_{j-1}}$ (cf. \eqref{lowerLIL}), a.s. for all large $j$, $\M_{n_{j-1}} \ge {1\over3} \log n_{j-1} \sim {\log 2\over3} j $, hence $x\equiv V(u) \gg  \lambda_j$ since $\lambda_j \sim a \log \log j$. Applying Lemma \ref{L:d1} gives that  on $\{ \M_{n_{j-1}} \ge {1\over3} \log n_{j-1}\}$, for some constant $C>0$, for all  $ \vert u \vert = n_{j-1}$, $$  \p\Big(  \M_{n_j - n_{j-1}} <  {3\over2} \log n_j + \lambda_j - x \Big)\big\vert_{x=V(u)} \le C \, V(u)   \, \ee^{- (V(u)- \lambda_j)}.$$

\noindent It follows that  \begin{eqnarray*}  \p\Big( A_j  \, \vert\, \F_{n_{j-1}}\Big) &\ge &  1_{(\M_{n_{j-1}} \ge {1\over3} \log n_{j-1})} \, \prod_{  \vert u \vert = n_{j-1}}  \Big( 1-   C   V(u) \, \ee^{- (V(u)- \lambda_j)} \Big)  \\
 	&\ge&   1_{(\M_{n_{j-1}} \ge {1\over3} \log n_{j-1})} \,   \exp\Big( - 2 C\, \sum_{ \vert u \vert = n_{j-1}} V(u) \, \ee^{- ( V(u) - \lambda_j)} \Big) \\
	&=&  1_{(\M_{n_{j-1}} \ge {1\over3} \log n_{j-1})} \,  \exp\Big( - 2 C \, \ee^{\lambda_j}  D_{n_{j-1}}  \Big) .
	\end{eqnarray*}

	Since   $D_{n_{j-1}} \to D_\infty$, a.s.,  and  $\ee^{\lambda_j} \sim (\log j)^a$ with $a<1$ ,   we get that  almost surely, $$ \sum_j  \p\Big( A_j  \, \vert\, \F_{n_{j-1}}\Big)= \infty,$$
	
	\noindent which according to Lévy's conditional form of Borel-Cantelli's lemma (\cite{levy}, Corollary 68), implies that $\p( A_i, \, \mbox{i.o.})=1$. Then  $$ \limsup_{n \to \infty} {1\over \log \log \log n} ( \M_n - {3\over2} \log n) \ge a, \qquad \mbox{ a.s.}$$  The lower bound follows by letting $a \to 1$. $\Box$

 \medskip
  {\noindent \it Proof of the upper  bound in Theorem \ref{t:3log}.} Let $\delta>0$ be small.  Recall \eqref{La}.   Let $a \ge a_0$ be   as in Lemma \ref{L:Nerman} such that  a.s. on $\S$, $\#\L_\lambda^{(a)} \ge \ee^{(1-\delta) \lambda}$ for all large $\lambda$.    Let $b>0$ such that $\ee^{-b} > q\equiv \p(\S^c)$.  By Lemma \ref{tightness}, there exists some constant $x_0>0$ such that $$ \p\Big( \max_{ n \le k \le 4 n} \M_k > { 3\over 2} \log n + x_0\Big) \le    \ee^{-b}, \qquad \forall n \ge n_0.$$

  \noindent Let $x_1:= x_0+a$. 
  Consider large integer $j$ and define $n_j:= 2^j$, $\lambda_j:=(1+ 2 \delta) \log \log\log n_j$. Define $$   B_j:= \Big\{ \max_{n_{j} < k \le n_{j+1}} \M_k > {3\over2} \log n_j+   \lambda_j + x_1 \Big\}\cap \S.$$

  Then,  \begin{eqnarray*} && \p\Big( B_j,  \#\L^{(a)}_{\lambda_j} \ge   \ee^{ (1-\delta) \lambda_j}, \max_{u \in \L^{(a)}_{\lambda_j}} \vert u \vert \le n_{j-1}  \Big)  \\
   &\le& \p\Big( \forall u \in \L^{(a)}_{\lambda_j}: \max_{n_{j-1} \le k \le n_{j+1}} \M_k^{(u)} >  {3\over2} \log n_j+  x_0,  \, \#\L^{(a)}_{\lambda_j} \ge   \ee^{(1-\delta)\lambda_j}\Big) 
   \\ &\le& \exp\Big( - b\,   \ee^{(1-\delta)\lambda_j}\Big) ,
  \end{eqnarray*}
  
  \noindent  whose sum on $j$ converges [$\delta$ being small].  On the other hand, by \eqref{infsup2}, $\p\big(\max_{u \in \L^{(a)}_{\lambda_j}} \vert u \vert    > n_{j-1}\big) \le  c_6 \ee^{- c_5 n_{j-1}^{1/3}} $ whose sum  again converges. Therefore, $ \sum_j \p \big( B_j,  \#\L^{(a)}_{\lambda_j} \ge   \ee^{ (1-\delta) \lambda_j}\big) < \infty$. By Borel-Cantelli's lemma,  almost surely, for all large $j$, the event $\{B_j,  \#\L^{(a)}_{\lambda_j} \ge  \ee^{(1-\delta)\lambda_j} \}$ does not hold; but we have chosen $a$ such that  on $\S$, $\#\L^{(a)}_{\lambda_j} \ge  \ee^{(1-\delta)\lambda_j} $ for all large $j$. Hence a.s. on $\S$, for all large $j$,    $  \max_{n_{j} < k \le n_{j+1}} \M_k \le  {3\over2} \log n_j+   \lambda_j + x_1, $ from which we  get that a.s. on $\S$, $$ \limsup_{ n \to \infty} {1\over \log\log\log n} ( \M_n - {3\over2} \log n) \le 1+2 \delta,$$
  
  \noindent yielding the upper bound as $\delta>0$ can be arbitrarily small. $\Box$

 \section{Proof of Theorem \ref{T:DEV1}} \label{S:5}
 

 \subsection{The B\"{o}ttcher case:  Proof of \eqref{dev2-main}}

 Recall \eqref{defll} for the stopping line  $\L_\lambda$.   
 
 \begin{lemma}[The B\"{o}ttcher case] \label{L:reduce} Under the same assumptions  as in Theorem \ref{T:DEV1}, for  any constant $a>0$, we have  \begin{equation} \label{numberL}
 \e \Big[ \ee^{ - a \# \L_\lambda}  \Big] =   \ee^{- \ee^{ (\beta+o(1)) \lambda}},  \qquad \lambda\to \infty. 
 \end{equation} \end{lemma}

 {\noindent\bf Proof of Lemma \ref{numberL}.}     Let us  check at first   the lower  bound in \eqref{numberL}.    Observe that $\p$-almost surely,  \begin{equation}\label{eq-d} D_\infty 
 =  \sum_{u \in \L_\lambda} \ee^{- V(u)} D_\infty(u),\end{equation}
 
 \noindent where conditioned on $\{V(u), u \in \L_\lambda\}$, $D_\infty(u)$ are independent copies of $D_\infty$.   Take $K_0$ large enough such that $ \e [ \ee^{- K_0 D_\infty} ] \le \ee^{- a}$, that is possible because $D_\infty >0$, $\p$-a.s. 
 Let $x= K_0\, \ee^{\lambda +K}$, where $K=\mbox{esssup} \max_{\vert u \vert =1} V(u) <\infty$ is as in  \eqref{bounded}.  Therefore \begin{eqnarray*}
\e \Big[ \ee^{- x D_\infty}\Big] & = & \e \Big[ \prod_{ u \in \L_\lambda} \e \big[ \ee^{- x \ee^{-y} D_\infty}\big]\big\vert_{y= V(u)\le \lambda+K}\Big]  \\
	&\le& \e \Big[ \prod_{ u \in \L_\lambda} \ee^{-a}\Big] = \e \Big[ \ee^{- a \#\L_\lambda}\Big].\end{eqnarray*}
 
\noindent  Hence $ \e \Big[ \ee^{- a \#\L_\lambda}\Big] \ge \e \Big[ \ee^{- x D_\infty}\Big] = \ee^{- x^{\beta+o(1)}}= \ee^{- \ee^{ (\beta+o(1)) \lambda}}$ gives the lower bound of \eqref{numberL}. 
 
 For the upper bound  of \eqref{numberL}, we use again \eqref{eq-d} to see that $D_\infty \le \ee^{-\lambda} \sum_{u \in \L_\lambda} D_\infty(u)$.  Take a constant $b>0$ such that $\e \big[ \ee^{- b D_\infty}\big] \ge \ee^{-a}$.  It follows that $$ \e \Big[ \ee^{- b \, \ee^\lambda\, D_\infty}\Big] \ge \e \Big[ \ee^{- b \sum_{ u \in \L_\lambda} D_\infty(u)}\Big] \ge \e \Big[ \ee^{-a \# \L_\lambda}\Big],$$ since conditioned on $\L_\lambda$, $(D_\infty(u))_{u \in \L_\lambda}$ are i.i.d. copies of $D_\infty$. Then \eqref{smalldev2} implies the upper bound of  \eqref{numberL}. $\Box$

 \medskip
 {\noindent \bf Proof of   \eqref{dev2-main}.}  By   Lemmas \ref{L:d1} and  \ref{tightness},   we can choose two positive  constants  $c_{11}$  and $c_{12}$ such that  for any $n \ge 1$,  \begin{eqnarray}  \min_{{n\over2} \le j \le n} \p \Big( \M_j \ge {3\over2}\log n - c_{11} \Big)  &\ge& \ee^{-c_{12}},  \label{tight-low}
 	\\   \p \Big( \max_{{n\over2} \le j \le  3 n} \M_j \ge {3\over2}\log n + c_{11} \Big)  &\le& \ee^{-c_{12}},  \label{tight-up}
	\end{eqnarray}
	
  For any $u \in \T$, define  as before $\M_j^{(u)}:= \min_{v \in \T_u, \vert v\vert= \vert u \vert +j} ( V(v)- V(u))$  for any $j \ge 0$. It follows  
  that \begin{eqnarray*}\p\Big( \M_n > { 3\over2} \log n + \lambda -c_{11} \Big) & \ge&  \p \Big( \forall u \in \L_\lambda,   \vert u \vert \le { n \over2}, M^{(u)}_{ n - \vert u \vert }  \ge  { 3\over2} \log n -c_{11}  \Big)  
	\\ &\ge& \e \Big[ \ee^{- c_{12} \, \# \L_\lambda} 1_{(     \max_{ u \in \L_\lambda} \vert u \vert \le { n \over2})}\Big] 
	\\&\ge&  \e \Big[ \ee^{- c_{12} \, \# \L_\lambda}  \Big] - \p\Big( \max_{ u \in \L_\lambda} \vert u \vert >{ n \over2}\Big)
	\\ &\ge &  \ee^{- \ee^{(\beta+o(1)) \lambda}} - c_6 \ee^{- c_5 \, n^{1/3}},
 \end{eqnarray*}

 \noindent by Lemma \ref{L:reduce} and \eqref{infsup2}.  The    lower bound in  \eqref{dev2-main} follows from the assumption that $\lambda=o(\log n)$.   
 
 To get  the upper bound in \eqref{dev2-main},   we  use the hypothesis \eqref{bounded} and obtain    that \begin{eqnarray*} && \p\Big( \max_{ n \le k \le 2n } \M_k> { 3\over2} \log n + \lambda +c_{11} +K \Big) \\
  & \le&  \p \Big( \forall u \in \L_\lambda, \max_{ u \in \L_\lambda} \vert u \vert \le { n \over2}, \max_{ n \le k \le 2n } \M^{(u)}_{ k - \vert u \vert } \ge  { 3\over2} \log n +c_{11}\Big)   +  \p\Big( \max_{ u \in \L_\lambda} \vert u \vert >{ n \over2}\Big)
  	\\ &\le& \e \Big[ \ee^{- c_{12} \, \# \L_\lambda}  \Big]  + c_6 \ee^{- c_5 \, n^{1/3}}  ,
 \end{eqnarray*}
 
 \noindent by \eqref{tight-up} and \eqref{infsup2}.    The upper bound follows from   Lemma \ref{L:reduce}.  
 $\Box$

\subsection{The Schr\"{o}der case: Proof of \eqref{dev1-main}}

   \medskip
 In the case   $q:=\p(\S^c)>0$, we need to estimate the probability that  the extinction happens after $\L_\lambda$:  
 
 \begin{lemma} \label{L:qexp} Assume \eqref{hyp}, \eqref{int1} and \eqref{defgamma1}.   Then for any $\lambda>0$,  $$ \p\Big( \{ \L_\lambda \neq \emptyset\}\cap \S^c\Big) = \e \Big[ q^{ \#\L_\lambda}\, 1_{( \#\L_\lambda>0)}\Big] \le q \,\ee^{ -\gamma \lambda}.$$ 
 \end{lemma}

 \medskip 
  
 {\noindent\bf Proof of Lemma \ref{L:qexp}.}  The above   equality is an immediate consequence of the branching property at the optional line $\L_\lambda$ (cf. \cite{BK04}).   
 
 To show the above inequality, we recall that  $\nu(u)$,  for any $u \in \T$,   denotes  the number of children of  $u$. Write $u < \L_\lambda$ if there exists some particle $v \in \L_\lambda$ such that $u < v$ [i.e. $u$ is an ancestor of $v$].  Then for the   tree up to $\L_\lambda$, the following equality holds: almost surely, \begin{equation}\label{tree-eq} \#\L_\lambda= 1 + \sum_{ \varnothing \le u < \L_\lambda} \big( \nu(u)- 1\big).\end{equation}

 \noindent Recall  \eqref{defgamma1}. Define a process $$ X_n:= \sum_{ \vert u \vert =n} \, \prod_{i=0}^{n-1} \big( q^{\nu(u_i)-1}\, 1_{(\nu(u_i)\ge1)}\big)\, \ee^{ \gamma V(u)}, \qquad n\ge 1,$$

 \noindent where as before, $u_i$ denotes the ancestor of $u$ at $i$th generation.  It is straightforward to check, by using the branching property, that $(X_n)_{n\ge 1}$ is a  (nonnegative) martingale with mean $1$.  Define $$X_{\L_\lambda}:= \sum_{ u \in \L_\lambda}\, \prod_{i=0}^{\vert u \vert -1}   \big( q^{\nu(u_i)-1}\, 1_{(\nu(u_i)\ge1)}\big)\, \ee^{ \gamma V(u)}, \qquad \lambda>0.$$ 
 
 According to Biggins and Kyprianou (\cite{BK04}, Lemma 14.1), $\e \big[ X_{\L_\lambda}\big]$ equals $\e[X_1]$ times  some probability term, hence $\e \big[ X_{\L_\lambda}\big] \le \e [X_1]=1$. 

 Notice that   for  any $u \in \L_\lambda$,  $\nu(u_i)\ge1$ for all $i < \vert u \vert$ and $\prod_{i=0}^{ \vert u \vert -1} \big( q^{\nu(u_i)-1}\, 1_{(\nu(u_i)\ge1)}\big)  = q^{ \sum_{ 0 \le i < \vert u \vert } ( \nu(u_i)- 1)} \ge q^{ \#\L_\lambda-1} $ by \eqref{tree-eq} [recalling $q<1$]. Then $ X_{\L_\lambda} \ge  q^{ \# \L_\lambda-1} \ee^{\gamma \lambda}$ on $\{\#\L_\lambda >0\}$. The Lemma follows from $\e[X_{\L_\lambda}]\le 1$. 
 $\Box$

\begin{lemma}\label{lowerproba=1} Assume \eqref{hyp}, \eqref{int1}, \eqref{defgamma1} and  \eqref{hyp-exp}. For any   $\delta>0$, there exist  an integer $m_\delta\ge 1$ and a constant $\lambda_0(\delta)>0$ such that   for all $\lambda\ge \lambda_0(\delta)$,   $$ \p\Big( 0 < \#\L_\lambda \le m_\delta \Big) \ge \ee^{ - ( \gamma+\delta ) \lambda}  . $$  \end{lemma}

{\noindent \bf Proof of Lemma \ref{lowerproba=1}:}   We discuss  the case $q=0$ and the case $q>0$ separately.

{\it (i)  First case: $q=0$.} We shall prove that \begin{equation}\label{low:q=0} \p \Big( \#\L_\lambda= 1\Big) \ge \ee^{- (\gamma+o(1)) \lambda}, \end{equation}

\noindent where as usual $o(1)$ denotes a quantity which goes to $0$ as $\lambda \to \infty$. 
To this end, we have by the change of measure (see Section \ref{S:2.2} and \eqref{WL}) that      \begin{equation} 
\p\Big( \#\L_\lambda=1 \Big)   =  \q \Big[ {1 \over W_{\L_\lambda}} 1_{( \#\L_\lambda=1    )}\Big] = \q \Big[ \ee^{  V( \w_{\tau_\lambda(\w)})}  1_{( \#\L_\lambda=1    )}\Big]   \ge   \ee^\lambda\, \q \Big( \#\L_\lambda=1   \Big).  \label{prob=1}
\end{equation}

\noindent  Notice that under $\q$, $\{  \#\L_\lambda=1\}$ means that $\L_\lambda= \{ \w_{ \tau_\lambda(\w)}\}$.   
Recall that  $\nu(u)$ denotes   the number of children of $u\in \T$.  Then   $  \q \Big( \#\L_\lambda=1    \big\vert \G_\infty\Big) = 1_{ (0\le  k < \tau_\lambda(\w),  \nu( \w_k)=1 )}  $  and thus \begin{equation} \label{lowpro1} \p\Big( \#\L_\lambda=1  \Big)   \ge   \ee^\lambda \, \q  \Big( 0\le  k < \tau_\lambda(\w),  \nu( \w_k)=1  \Big)   .
	\end{equation}

Recall  \eqref{defgamma2} for $\gamma$. We claim that   \begin{equation} \label{qproba1} \q  \Big( 0\le  k < \tau_\lambda(\w),  \nu( \w_k)=1 \Big) = \ee^{ - (1+\gamma +o(1)) \lambda}. \end{equation}

\noindent 
To get  \eqref{qproba1}, we  use the fact   (cf. Section \ref{S:2.2}) that  $(   \sum_{ u \in \mho( \w_k)} \delta_{\{ \Delta V(u)\}},  \Delta V(\w_k))_{k \ge 1}$ are i.i.d. under $\q$, where $\Delta V(u):= V(u) - V({\buildrel \leftarrow \over u}) $ for any $ u \not= \varnothing \equiv \w_0$.    Notice that $\nu(\w_{k-1})=  1 + \#   \mho( \w_k)  $.

Let us check that   the process   $$ U_n:=  \ee^{ (1+\gamma) V(\w_n)} 1_{( \forall  1\le k \le   n, \nu(\w_{k-1})=1)}, \qquad n \ge 1,$$  

\noindent is a $\q$-martingale of mean $1$.  In fact, $U_n$ is a product of $n$ i.i.d.  variables, then  it is enough to check that $  \q\big[U_1\big]=1$. But $\q\big[U_1\big]=   \q \big[ \ee^{ (1+\gamma) V(\w_1)} 1_{( \nu(\w_0)=1)}\big] = \e \big[ \sum_{ \vert u \vert =1} \ee^{  \gamma  V(u) } 1_{ (\nu=1)} \big]  =1, $ as claimed. By the optional stopping  theorem and the Fatou lemma, we get that $ \q\big[ U_{ \tau_\lambda(\w)}\big] \le 1$, which implies the upper bound in \eqref{qproba1} since $ V( \tau_\lambda(\w)) > \lambda$  [under $\q$, $\tau_\lambda(\w)$ is a.s. finite].

To get the lower bound in \eqref{qproba1}, let $\varepsilon>0$ be small. Fix some large constant $C $ whose value will be determined later.       Let us find some $\gamma_C $ such that the process $$U^{(C)}_n:= \ee^{ (1+\gamma_C) V(\w_n)} 1_{( \forall  1\le k \le   n, \nu(\w_{k-1})=1 , \,  \Delta V(\w_k) \le C)}, \qquad n \ge 1,$$  

\noindent is a $\q$-martingale with mean $1$. As for $U_n$, the constant  $\gamma_C $\footnote{For the existence of such constant,  we used the integrability assumption \eqref{hyp-exp}: the convex function $f: b \to \e \big[ \sum_{\vert u \vert =1} \ee^{b V(u)}  1_{(\nu=1)} \big] $  has a derivative $f'(\gamma) \ge {f(\gamma)- f(0) \over \gamma} >0$ hence $f$ is increasing at $\gamma$. Then $f(a)>f(\gamma)=1$. Take $C_0$ large enough such that $\e \big[ \sum_{\vert u \vert =1} \ee^{a V(u)} 1_{( \nu=1 , V(u) \le C_0)} \big]>1$, then   such $\gamma_C$ exists for all $C\ge C_0$. We shall use the existences   of similar constants  later without further explanations.}   is determined by $$ 1=   \e \Big[ \sum_{ \vert u \vert =1} \ee^{  \gamma_C  V(u) } 1_{ (\nu=1,  V(u) \le C)} \Big] ,$$ 

\noindent  where  for $|u|=1$, $\Delta V(u)= V(u)$.  Plainly $\gamma_C \to \gamma$ as $C \to \infty$.  Choose $C$ sufficiently large such that $\gamma_C \le \gamma+\varepsilon  $. Since $(U^{(K)}_k, k \le \tau_\lambda(\w))$ is uniformly bounded by $\ee^{(1+\gamma_C) (\lambda+C)}$. By the optional stopping  theorem, we obtain that $$ 1= \q \Big[ U^{(C)}_{\tau_\lambda(\w)}\Big]  \le 
\ee^{ (1+\gamma_C) (\lambda+C)}  \q\Big( \forall 1\le k \le   n, \nu(\w_{k-1})=1\Big),$$ finishing the proof of \eqref{qproba1} as $\varepsilon$ can be arbitrarily small.  The  Lemma (in the case $q=0$) follows  from \eqref{qproba1}  and \eqref{lowpro1}.

{\it (ii) Second (and last) case: $q>0$.}   We can not repeat the same proof as before, for instance  $p_1\equiv \p(\nu=1)$ may vanish. 

Again  by the change of measure we have that for any integer $m\ge1$, \begin{equation} \label{prob=m}
\p \Big( 0 < \# \L_\lambda \le m \Big) = \q \Big[ { 1 \over W_{\L_\lambda}} 1_{(  \#\L_\lambda \le m)} \Big] \ge {1 \over m} \ee^\lambda\, \q \Big(   \# \L_\lambda \le m\Big), \end{equation}

\noindent where we used the facts  that $W_{\L_\lambda} = \sum_{ u \in \L_\lambda} \ee^{- V(u)} \le m \ee^{-\lambda}$ on $\{ \#\L_\lambda\le m\}$ and under $\q$, $\L_\lambda $ contains at least the singleton $\{\w_{\tau_\lambda(\w)}\}$.  Define for any $x>0$,  $$q(x):= \p \big( \sup_{v \in \T} V(v) \le x\big) = \p \big( \L_x= \emptyset\big),  $$

\noindent    with the usual  convention that $\sup_\emptyset = 0$.  Plainly,  $ \lim_{ x \to \infty} q(x)= \p \big( \sup_{v \in \T} V(v) < \infty \big)=  \p \big( \S^c\big) = q.$   For any small $\varepsilon>0$, there exists some $x_0=x_0(\varepsilon)>0$ such that $ q(x) \ge q- \varepsilon$ for all $x \ge x_0$. 



Let $\delta>0$ be small.  Before bounding below $\q \big(   \# \L_\lambda \le m\big)$ with some $m=m_\delta$, we first choose some constants.  Let $\alpha$ be   large and $\varepsilon$ be small   whose values  will be determined  later. Recall that $\mho( \w_k)$ denotes the set of brothers of $\w_k$.  Let us    choose   a constant $\gamma_{\alpha,   \varepsilon} $ such that \begin{equation} \label{uae} U^{(\alpha,  \varepsilon)}_n:= \ee^{(1+ \gamma_{\alpha,  \varepsilon} )\, V( \w_n)}\, (q-\varepsilon)^{\sum_{ 0 \le k < n} ( \nu(\w_k)- 1)} 1_{( \forall k < n, \forall u \in \mho( \w_k), \, \Delta V(u) \le \alpha )}, \qquad n \ge1, \end{equation}

\noindent is a $\q$-martingale with mean $1$.  As before, such $\gamma_{\alpha,  \varepsilon}$ is determined by the following equalities    \begin{eqnarray*}  1 &= & \q\Big[ \ee^{ (1+\gamma_{\alpha,   \varepsilon}) V(\w_1)} (q-\varepsilon)^{ \nu(\w_0)-1}\, 1_{( \max_{ \vert u \vert  =1, u \neq \w_1} V(u) \le \alpha)}\Big]  \\
	&=& \e \Big[ \sum_{ \vert u \vert=1} \ee^{ \gamma_{\alpha,  \varepsilon} V(u)}  (q-\varepsilon)^{ \nu -1}\, 1_{( \max_{ \vert v \vert  =1, v \neq u } V(v) \le \alpha)}\Big] .  \end{eqnarray*}


The existence of $\gamma_{a, \varepsilon}$ follows from  \eqref{defgamma1} and   \eqref{hyp-exp}. Clearly  $\gamma_{\alpha, \varepsilon} \to \gamma$ as $\alpha \to \infty$ and $\varepsilon \to 0$.  Fix now   $\alpha\equiv \alpha(\delta)>0  $ (large enough)    and   $\varepsilon\equiv \varepsilon(\delta) >0$ (small enough) such that  $\gamma_{\alpha, \varepsilon} < \gamma +\delta$.  Choose a constant   $x_0 \equiv x_0(\delta) >0$ such that $ q(x) \ge q- \varepsilon$ for all $x \ge x_0$.

On the other hand,  we remark that   \eqref{hyp} and  \eqref{defgamma1} imply   that \begin{equation}\label{newhyp}   \p\Big( 1\le \nu< \infty, \, \max_{\vert u \vert =1} V(u) >0\Big) >0. \end{equation} 

\noindent In fact, $ \e \big[ 1_{(1\le \nu < \infty)} q^{\nu-1} \sum_{\vert u \vert=1} \ee^{\gamma V(u)} 1_{( V(u)>0)}\big] = 1-  \e \big[ 1_{(1\le \nu < \infty)} q^{\nu-1} \sum_{\vert u \vert=1} \ee^{\gamma V(u)} 1_{( V(u)\le 0)}\big] >  1-  \e \big[ 1_{(1\le \nu < \infty)} q^{\nu-1}  \nu \big]  > 0$, hence \eqref{newhyp} holds.  It follows that  there are    some integer $n_*\ge 1$ and some positive constants $c_* $ and $b_*$ such that  \begin{equation}\label{n*c*} b_*\le \e \Big[ 1_{( \nu \le n_*)} \, \sum_{\vert u \vert =1} \ee^{- V(u)} 1_{( V(u) \ge c_*)}\Big] = \q \Big( \nu(\w_0)\le n_*, \,   V(\w_1) \ge c_*\Big),
\end{equation} where the last equality follows from the change of measure formula (Section \ref{S:2.2} (i), $\w_0= \varnothing$).  

Choose (and  fix) a constant $L\ge   \alpha +x_0$ such that ${L\over c_*}$ is an integer.  Define $m_\delta:= (n_*)^{L/c_*}$. 
   Recall \eqref{tau} for the definition of $\tau_\lambda(u)$. For  any $\lambda> 2L$, we consider the following events  \begin{eqnarray*}  A_1&:=&  \Big \{   \forall k < \tau_{\lambda-L}(\w), \, \forall u \in \mho( \w_k),   \Delta V(u)  \le \alpha,  \L_\lambda^{(u)}= \emptyset \Big\}  ,    \\
 A_2 &:=& \Big \{   \forall   \tau_{\lambda-L}(\w) \le k <  \tau_{\lambda-L}(\w) + { L\over c_*},  \, \forall u \in \mho( \w_k),   \nu(u)=0, \, \nu(\w_{k-1}) \le n_*, \Delta V(\w_k) \ge c_*  \Big\} ,   \end{eqnarray*} 

\noindent  where  $\L_\lambda^{(u)}:=  \T_u \cap \L_\lambda$ and  $\nu(u)$ denotes the number of children of $u$.

Observe that on $A_1\cap A_2$,    $\tau_\lambda(\w) \le  \tau_{\lambda-L}(\w) + { L\over c_*}$, and  $\#\L_\lambda   \le  (n_*)^{L/c_*}\equiv m_\delta$.  Since $q>0$, $p_0\equiv \p(\nu=0) >0$, it follows from the spinal decomposition (Section \ref{S:2.2} (iii))  that   
\begin{eqnarray}   \q \Big( \# \L_\lambda \le m_\delta\Big) & \ge&  \q \Big( A_1 \cap A_2 \Big)  \nonumber \\ 	&=& \q \Big[  B_1   \, \prod_ {k=\tau_{\lambda-L}(\w) }^{  \tau_{\lambda-L}(\w) + { L\over c_*}-1}   \prod_{u \in \mho( \w_k) } p_0\times 1_{( \nu(\w_{k-1}) \le n_*, \Delta V(\w_k) \ge c_*  )} \Big ] \nonumber
	\\ & \ge& p_0^{m_\delta}\, \q \Big[  B_1   \, \prod_ {k=\tau_{\lambda-L}(\w) }^{  \tau_{\lambda-L}(\w) + { L\over c_*}-1}   1_{( \nu(\w_{k-1}) \le n_*, \Delta V(\w_k) \ge c_*  )} \Big ] ,  \label{lowqproba3}  \end{eqnarray}

\noindent where $$ B_1:= \prod_{k < \tau_{\lambda-L}(\w)} \prod_{u \in \mho( \w_k)} q(\lambda-V(u)) 1_{(\Delta V(u) \le \alpha)}  \ge \prod_{k < \tau_{\lambda-L}(\w)} (q-\varepsilon)^{ \nu(\w_{k-1})- 1}  1_{(\max_{ u \in \mho( \w_k)} \Delta V(u) \le \alpha)} =: B_2 ,$$

 \noindent by using the fact that  for any $u \in \mho( \w_k)$ with $k< \tau_{\lambda-L}(\w)$, $V(u) \le \lambda-L +\alpha \le \lambda - x_0$, and   $q(\lambda-V(u)) \ge q(x_0) \ge q-\varepsilon$. 
 
 Recall that under $\q$,  $(   \sum_{ u \in \mho( \w_k)} \delta_{\{ \Delta V(u)\}},  \Delta V(\w_k))_{k \ge 1}$ are i.i.d.; then the strong Markov property implies that under $\q$ and conditioned on $\G_{\tau_{\lambda-L}(\w)}$, $(\nu(\w_{k-1}), \Delta V(\w_k))_{k \ge  \tau_{\lambda-L}(\w)}$  are i.i.d.,  of common law that of $(\nu(\w_0), V(\w_1))$.  Therefore,   \begin{equation}\label{lowqproba4} \q \Big( \# \L_\lambda \le m_\delta\Big)   \ge   p_0^{m_\delta}\, \q \Big[  B_2 \Big]\,    \, \q \Big( \nu(\w_0) \le n_*,   V(\w_1) \ge c_*  \Big)^{L/c_*}  \ge p_0^{m_\delta}\, b_*^{L/c_*}\, \q \Big[  B_2 \Big].    \end{equation}

It remains to estimate $\q[B_2]$.  Going back to \eqref{uae} and  applying   the optional stopping  theorem at $\tau_{\lambda-L}$ for $U^{(\alpha, \varepsilon)}$ (which remains bounded up to $\tau_{\lambda-L}$), we get   that $$   \q [B_2]=  \q \Big[  (q-\varepsilon)^{ \sum_ { 0 \le k < \tau_{ \lambda-L}(\w)} (\nu( \w_k)-1)}     1_{(   \forall k < n, \forall u \in \mho( \w_k), \, \Delta V(u) \le \alpha )  }  \Big]  \ge  \ee^{- (1+\gamma_{\alpha, \varepsilon}) (\lambda- L  +\alpha)}.$$

  In view of  \eqref{prob=m} and \eqref{lowqproba4}, this implies that $$\p \Big( 0< \# \L_\lambda \le m_\delta\Big) \ge  {1\over m_\delta}\, p_0^{m_\delta}\, b_*^{L/c_*}\, \ee^{L-\alpha}\, \ee^{-  \gamma_{\alpha, \varepsilon} (\lambda-L+\alpha)  }.$$

\noindent Then we have proved  the Lemma in the case $q>0$ [by   choosing  a sufficiently large $\lambda_0(\delta)$]. $\Box$

 \begin{lemma}[The Schr\"{o}der case] \label{L:numb2} Under the same assumptions  as in Theorem \ref{T:DEV1}, for  any constant $a>0$, we have  \begin{equation} \label{numb2}
 \e \Big[ \ee^{ - a\, \# \L_\lambda} 1_{( \#\L_\lambda>0)}\Big] =     \ee^{ - (\gamma+o(1)) \lambda},  \qquad \lambda\to \infty. 
 \end{equation} \end{lemma}

{\noindent\bf Proof of Lemma \ref{L:numb2}.} From  Lemma \ref{lowerproba=1}, the lower bound of  \eqref{numb2} follows immediately.
We also mention that in the cases when  $q=0$ or $q>0$ but    $ 0<a < \log (1/q)$, we can give a proof of the lower bound of \eqref{numb2} in the same way as that of \eqref{numberL}.

 For the upper bound, we proceed in the same way as  in the proof of Lemma \ref{L:reduce}, but by paying attention to the possibility of extinction of the system. Take $b>0$ such that  $\e \big[ \ee^{- b D_\infty}\big] \ge \ee^{-a}$. By \eqref{eq-d},  $\ee^\lambda\, D_\infty \le \sum_{u \in \L_\lambda} D_\infty(u)$, then \begin{eqnarray*}  \e \Big[ \ee^{- b \, \ee^\lambda\, D_\infty} 1_{(D_\infty >0)} \Big] &\ge & \e \Big[ \ee^{- b \sum_{u \in \L_\lambda} D_\infty(u)} 1_{(D_\infty >0)} \Big] 
 	\\& \ge&  \e \Big[ \ee^{- b \sum_{u \in \L_\lambda} D_\infty(u)} 1_{( \#\L_\lambda>0)} \Big]  - \p\Big( \{\#\L_\lambda>0\} \cap \S^c\Big) 
 	\\& \ge&    \e \Big[ \ee^{ - a\, \# \L_\lambda} 1_{( \#\L_\lambda>0)}\Big]  - \p\Big( \{\#\L_\lambda>0\} \cap \S^c\Big) .
	\end{eqnarray*}

\noindent By \eqref{smalldev1new}, $ \e \big[ \ee^{- b \, \ee^\lambda\, D_\infty} 1_{(D_\infty >0)} \big]  \le C \ee^{- \gamma \lambda}$,  which together with Lemma \ref{L:qexp} yield the upper bound in \eqref{numb2}. $\Box$

 \medskip

We now are ready to give the proof of  \eqref{dev1-main}:

\medskip
{\noindent\bf Proof of   \eqref{dev1-main}.}   Let us prove at first the the lower bound in \eqref{dev1-main}.  By Lemma \ref{L:d1}, there are $c_{13}>0$ (large enough) and $c_{14}>0$ (small enough) such that $ \min_{ {n\over2} \le k \le n} \p ( \M_k \ge {3\over2} \log n - c_{13}, \S) \ge c_{14}$ for all $n \ge 1$.

Let $\delta>0$ be small and let $m_\delta\ge1$ and $\lambda_0(\delta)>0$ be as in Lemma \ref{lowerproba=1}.  Let $\lambda \ge \lambda_0(\delta)$. Remark that  $$ \p \Big( \M_n > { 3\over2} \log n +   \lambda- c_{13},  \S \Big)   \ge \p\Big( 0< \#\L_\lambda \le m_\delta, \forall u \in \L_\lambda,  \M^{(u )}_{ n-\vert u \vert} >  { 3\over2} \log n - c_{13}, \vert u \vert \le {n \over 2}  , \,  \S^{(u )}\Big), $$

\noindent where as before, $S^{(u)}= \{ \T_u \mbox{ suvives}\}$ and $\M_j^{(u)}:= \min_{v\in \T_u, \vert v \vert=\vert u \vert +j} (V(v)- V(u))$ for any $j \ge 0$.     It follows that \begin{eqnarray*} \p \Big( \M_n > { 3\over2} \log n +   \lambda -c_{13},  \S \Big)   &\ge  & \, (c_{14})^{m_\delta}\, \p\Big( 0< \#\L_\lambda \le m_\delta, \max_{u \in \L_\lambda} \vert u \vert  \le {n \over 2} \Big)   	\\& \ge& (c_{14})^{m_\delta}\,  \Big(  \p\Big(0  < \#\L_\lambda \le m_\delta\Big) -\p\Big(  \max_{u \in \L_\lambda} \vert u \vert  >  {n \over 2} \Big) \Big) 
	\\ &\ge& (c_{14})^{m_\delta}\,  \Big( \ee^{- (\gamma+\delta) \lambda} - c_6 \ee^{- c_5 n^{1/3}}\Big), 
	\end{eqnarray*} by Lemma \ref{lowerproba=1} and  \eqref{infsup2}. The lower bound of \eqref{dev1-main} follows.  
	
	We prove now the upper bound in  \eqref{dev1-main}.  By assumption \eqref{hyp-exp} holds for any $a>0$, hence $S_1$ has all exponential moments. It follows from \eqref{overshoot}  that  for any $a>0$, there exists some $C_a>0$ such that  \begin{equation} \label{expmoment1}  
   \p \Big( S_{\tau_\lambda  }  - \lambda \ge x \Big) \le C_a\,  \ee^{- a x}, \qquad \forall \, x\ge 0. \end{equation}
  
 \noindent  Let $\delta>0$ be small and $a>(1+\gamma)/\delta +1.$   Then  \begin{equation} \p\Big( \max_{u \in \L_ \lambda} V(u) > (1+ \delta )  \lambda \Big)  \le \e \Big[ \sum_{  u \in \L_\lambda} 1_{( V(u) > (1 + \delta )  \lambda)} \Big] =\e\Big[ \ee^{ S_{\tau_\lambda}} 1_{( S_{\tau_\lambda \ge (1+\delta) \lambda)}}\Big] =o(\ee^{- \gamma \lambda}),\label{expmoment2}
  \end{equation}

  \noindent where the last equality follows easily from \eqref{expmoment1}.    Define \begin{equation}\label{uan1} A_{\eqref{uan1}} := \Big\{ \max_{u \in \L_ \lambda} V(u)  \le (1+ \delta )  \lambda, \, \max_{u \in \L_\lambda} \vert u \vert \le { n\over 2}\Big\}.
  \end{equation}
  
 \noindent    Then by \eqref{infsup2},  for all large $n\ge n_0$ and $0 < \lambda =o( \log n)$, \begin{eqnarray*} \p\Big( A_{\eqref{uan1}}^c\Big) &\le&  \p\Big( \max_{u \in \L_ \lambda} V(u) > (1+ \delta )  \lambda \Big) + \p\Big(   \max_{u \in \L_\lambda} \vert u \vert  > { n\over 2}\Big)  \\
 	&\le& o(\ee^{-\gamma \lambda}) + c_6\,  \ee^{- c_5 \, n^{1/3}}= o(\ee^{-\gamma \lambda})  .
 \end{eqnarray*}

   On $\S\cap \{ \M_n > { 3\over2} \log n + (1+2\delta)  \lambda\}$, $\L_\lambda \neq \emptyset $.   Consider $\lambda$ such that $\delta \lambda  < \log n$.    Therefore,  \begin{eqnarray}  && \p \Big(  \max_{ n \le k \le 2 n}  \M_k   > { 3\over2} \log n + (1+2\delta)  \lambda,   \, \S\Big)   \nonumber
  	\\ &\le & \p \Big(  \max_{ n \le k \le 2 n}  \M_k   > { 3\over2} \log n + (1+2\delta)  \lambda, A_{\eqref{uan1}}  ,  \L_\lambda\neq \emptyset \Big) + o(\ee^{-\gamma \lambda})   \nonumber \\
  	&\le&  \p \Big( \forall u \in \L_\lambda, \max_{{n\over 2}  \le j \le 2 n  } \M^{(u)}_j> { 3\over 2} \log n + \delta \lambda ,  \L_\lambda\neq \emptyset \Big) + o(\ee^{-\gamma \lambda})  \nonumber \\
	&=:& B_{\eqref{dev1}}+ o(\ee^{-\gamma \lambda})  , \label{dev1}
  \end{eqnarray}

  \noindent where $\M^{(u)}_k:= \max_{v \in \T_u, \vert v\vert= \vert u \vert + k} \big( V(v)- V(u)\big)$. Conditioning on  $\F_{\L_\lambda}$, $\M^{(u)}_\cdot$ are i.i.d. copies of $\M_\cdot$.  By Lemma \ref{tightness} (with $a=4$), there exist  some $c_{15}>0$ and $\lambda_0$ such that   ($\delta$ being fixed) for all large $n \ge n_0(\lambda_0)$, $$  \p\Big(  \max_{{n\over 2} \le k \le 2 n} \M_k \ge { 3\over 2} \log n + \delta \lambda_0 \Big) \le \p( \S^c) +  \p^*\Big(  \max_{{n\over 2} \le k \le 2 n} \M_k \ge { 3\over 2} \log n + \delta \lambda_0 \Big)  \le \ee^{ -c_{15}}.$$ 
  
  \noindent Then by conditioning on $\F_{\L_\lambda}$, we get that \begin{eqnarray*}  B_{\eqref{dev1}} &\le& \e \Big[ \ee^{- c_{15}\, \# \L_\lambda} 1_{(  \L_\lambda\neq \emptyset )} \Big]  = \ee^{- (\gamma+o(1)) \lambda},  \end{eqnarray*}

  \noindent by Lemma \ref{L:numb2}.   This  and \eqref{dev1} prove  the upper bound in \eqref{dev1-main} since $\delta$ can be arbitrarily small. $\Box$

         \medskip
    
{\noindent\bf Acknowledgments.} We are very grateful to two anonymous referees for their careful readings and helpful comments on the first version of this paper. We also thank  Zhan Shi for the reference \cite{G00}.


\begin{thebibliography}{13}
\baselineskip=10pt

 
 
 \bibitem{AbR}   Addario-Berry, L. and Reed, B. (2009).
 Minima in branching random walks.
  {\it Ann.  Probab.}  {\bf 37}  pp. 1044--1079.


\bibitem{A11} A{\"{\i}}d{\'e}kon, E. (2013).
     Convergence in law of the minimum of a branching random walk.      {\it Ann. Probab.}    {\bf 41} pp. 1362--1426.
 
 


 
\bibitem{EOY} A{\"{\i}}d{\'e}kon, E., Y. Hu and O. Zindy. (2013). The precise tail behavior of the total progeny of a killed branching random walk. {\it Ann. Probab.}   {\bf 41}   pp. 3786--3878.

 
\bibitem{AS12} A{\"{\i}}d{\'e}kon, E. and Shi, Z. (2014). The Seneta-Heyde scaling for the branching random walk. {\it Ann. Probab.}    {\bf  42}  pp. 959--993.



 
 \bibitem{AN} Athreya, K.B and Ney, P.E.  (1972). {\it Branching processes.} Springer-Verlag, Berlin, New-York.  
 



\bibitem{BGMS} 
  Berestycki, N.,    Gantert, N.,    Morters, P. and    Sidorova, N. (2014) Galton-Watson trees with vanishing martingale limit.  {\it J. Stat. Phys. } {\bf 155}  pp.  737--762.

\bibitem{B76}
    Biggins, J.D. (1976).
    The first- and last-birth problems for a
    multitype age-dependent branching process.
    {\it Adv. Appl. Probab.} {\bf 8}, 446--459.



\bibitem{BK04}
    Biggins, J.D. and Kyprianou, A.E. (2004).
    Measure change in multitype branching.
    {\it Adv. Appl. Probab.} {\bf 36}, 544--581.


 

 \bibitem{BZ06}
    Bramson, M.D. and Zeitouni, O. (2009).
    Tightness for a family of recursion equations. {\it Ann. Probab.,} {\bf 37},          615--653
    



\bibitem{Chang94}  Chang, J.T. (1994). Inequalities for the overshoot.  {\it   Annals  Appl. Probab.} {\bf  4}  No. 4   pp. 1223--1233.

\bibitem{CRW}
    Chauvin, B., Rouault, A. and Wakolbinger, A.  (1991).  Growing conditioned trees.  {\it Stoch. Proc. Appl.} {\bf 39}, 117--130.



 
 \bibitem{D80}
 Doney, R.A.  (1980).
 Moments of ladder heights in handom walks.
 {\it J. Appl. Probab.},   {\bf 17}, {248--252}.







\bibitem{FZ10}  Fang, M. and Zeitouni, O. (2010). Consistent minimal displacement of branching  random walks. {\it Electronic Comm. Probab.}   106--115.

\bibitem{FHS}
Faraud, G., Hu, Y. and Shi, Z. (2012).  Almost sure convergence for stochastically biased
random walks on trees. {\it Prob. Th. Rel. Fields.} {\bf 154} 621--660. 


 \bibitem{feller}
   Feller, W. (1971).
    {\it An Introduction to Probability Theory and
   its Applications.} Vol. II. Second edition.
     Wiley, New York.

\bibitem{FW07}   Fleischmann,  K. and   Wachtel, V. (2007). 
Lower deviation probabilities for supercritical Galton-Watson processes. {\it Ann. I. H. Poincaré} {\bf 43} 233--255.


\bibitem{FW09}   Fleischmann,  K. and   Wachtel, V. (2009).  On the left tail asymptotics for the limit law of supercritical Galton-Watson processes in the B\"{o}ttcher case.  {\it Ann. I. H. Poincaré} {\bf 45} 201--225. 

\bibitem{G00} Gatzouras, D. (2000).  On the lattice case of an almost sure renewal theorem for branching random walk. {\it Adv. Appl. Proba.} {\bf 32} 720--737.


\bibitem{H74}
    Hammersley, J.M. (1974).
    Postulates for subadditive processes.
    {\it Ann. Probab.} {\bf 2}, 652--680.

\bibitem{HR} Harris, S.C. and Roberts, M.I. (2011+).  The many-to-few lemma and multiple spines. (Preprint, arXiv:1106.4761)


  
\bibitem{H12}
Hu, Y. (2012+). The almost sure limits of the minimal position and the additive martingale
in a branching random walk. {\it J. Theor. Probab.} (to appear)  available at arXiv:1211.5309


\bibitem{HS09}
    Hu, Y.\ and Shi, Z.  (2009).
    Minimal position and critical martingale
    convergence in branching random walks,
    and directed polymers on disordered
    trees.
    {\it Ann.\ Probab.}, {\bf 37}, 742--789.

\bibitem{J12}
    Jaffuel, B.\ (2012).
    The critical barrier for the survival of
    the branching random walk with absorption.
  To appear in  {\it Annales de l'Institut Henri Poincar\'e.}   {\bf 48}  pp. 989--1009.    
  
 \bibitem{K75}
    Kingman, J.F.C. (1975).
    The first birth problem for an age-dependent
    branching process.
    {\it Ann. Probab.} {\bf 3}, 790--801.




\bibitem{levy} Lévy, P. (1937).  {\it Théorie de l'addition des variables aléatoires.}  Gauthier-Villars, Paris.




\bibitem{Liu99} Liu, Q.S. (1999). Asymptotic properties of supercritical age-dependent branching processes and homogeneous branching
random walks.  {\it Stoch. Proc. Appl.}  {\bf 82}  61--87.

 \bibitem{Liu01}   Liu, Q.S. (2001).   Asymptotic properties and absolute continuity of   laws stable by random weighted mean.  {\it Stoch. Proc. Appl.}  {\bf 95}, 83--107.


\bibitem{Lorden} Lorden, G. (1970). On excess over the boundary. {\it   Ann. Math. Stat.} {\bf 41}  520--527.

\bibitem{L97}
    Lyons, R. (1997).
    A simple path to Biggins' martingale convergence
    for branching random walk.
    In: {\it Classical and Modern Branching
    Processes} (Eds.: K.B.~Athreya and P.~Jagers).
    {\it IMA Volumes in Mathematics and its
    Applications} {\bf 84}, 217--221. Springer, New
    York.




\bibitem{Mallein} Mallein, B. (2016+).  Asymptotic of the maximal displacement in a branching random walk.  https://arxiv.org/abs/1605.08292

\bibitem{M74}   Mogulskii, A.~A.\ (1974).    Small deviations in the space of     trajectories.    {\it Theory Probab.\ Appl.} {\bf 19},    726--736.
 
 
 \bibitem{N81} Nerman, O. (1981).  On the convergence of supercritical general (C-M-J) branching processes. {\it Z. Wahrsche. verw. Gebiete}  {\bf 57}  365--395.
 
 
 
\bibitem{Shi} Shi, Z. (2015). {\it Branching Random Walks.}     Lecture notes from the 42nd Probability Summer School held in Saint Flour, 2012. Lecture Notes in Mathematics, 2151. École d'Été de Probabilités de Saint-Flour. Springer, Cham.

  
\end{thebibliography}
\end{document}